\documentclass[12pt]{article} 

\usepackage[all]{xy}
\usepackage{amssymb}
\usepackage{graphicx}
\usepackage{amsthm}

\newenvironment{pf}{\begin{proof}}{\end{proof}}

\newtheorem{thm}{Theorem}[section]
\newtheorem{lemma}[thm]{Lemma}
\newtheorem{cor}[thm]{Corollary}
\newtheorem{prop}[thm]{Proposition}

\newtheorem{defin}[thm]{Definition}

\newtheorem{example}[thm]{Example}
\newtheorem{con}[thm]{Construction}
\newtheorem{notn}[thm]{Notation}

\newenvironment{ex}{\begin{example}\rm}{\end{example}\rightline{$\circ$}}
\newenvironment{const}{\begin{con}\rm}{\end{con}}
\newenvironment{notation}{\begin{notn}\rm}{\end{notn}}

\def\C{{\mathbb C}}

\def\P{{\mathbb P}}

\def\R{{\mathbb R}}
\def\Z{{\mathbb Z}}

\def\cA{{\cal A}}
\def\cB{{\cal B}}

\def\E{{\cal E}}

\def\cI{{\cal I}}

\def\cL{{\cal L}}

\def\N{{\cal N}}
\def\cO{{\cal O}}

\def\T{{\cal T}}

\def\V{{\cal V}}

\def\cX{{\cal X}}
\def\cY{{\cal Y}}
\def\cZ{{\cal Z}}

\def\a{\alpha}
\def\b{\beta}
\def\suma{\sigma}
\def\product{\prod}
\def\L{\cL}
\def\d{\partial}

\def\hW{{\widehat W}}
\def\hM{{\widehat M}}
\def\hN{{\widehat N}}
\def\ratmap{- \rightarrow}

\def\a{\alpha}

\def\iso{\cong}

\def\lra{\longrightarrow}

\def\ra{\rightarrow}

\def\operatorname#1{\mathop{\rm #1}\nolimits}

\def\Proj{\operatorname{Proj}}

\def\Hom{\operatorname{Hom}}

\def\Spec{\operatorname{Spec}}

\def\im{\operatorname{im}}

\newcommand\fiveandhalfroottwo{5.71}
\newcommand\threeandhalfroottwo{3.71}
\newcommand\rootthree{1.73205}

\addtocounter{section}{-1}

\title{On phylogenetic trees -- a geometer's view \thanks{Deticated to
Andrzej Bia{\l}ynicki-Birula}}

\author{Weronika Buczy\'nska 
and Jaros\l{}aw
A. Wi\'sniewski\thanks{Partially supported by Polish KBN (grant 1P03A03027)
and MPIM, Bonn, Germany}\\
Instytut Matematyki UW, Banacha 2, 02-097 Warszawa, Poland }

\begin{document}

\maketitle

\noindent{\bf Abstract:}\ \ In the present note we investigate
projective varieties which are geometric models of binary symmetric
phylogenetic 3-valent trees. We prove that these varieties have
Gorenstein terminal singularities (with small resolution) and they are
Fano varieties of index 4. Moreover any two such varieties associated
to trees with the same number of leaves are deformation equivalent,
that is, they are in the same connected component of the Hilbert
scheme of the projective space. As an application we provide a simple
formula for computing their Hilbert-Ehrhard polynomial.

\tableofcontents

\section{Introduction}

Algebraic geometry, a classical, almost ancient, branch of pure
mathematics, is constantly stimulated by questions arising in
applicable mathematics and other sciences. String theory and mirror
conjecture from mathematical physics, coding theory and image
recognition from computer sciences --- to mention just a few of the
big areas of sciences which had enormous impact on the development of
algebraic geometry in the past decade. Now the modern biology with its
computational aspects and relations to statistics seems to be making
its way into this branch of mathematics.

Although the roots of questions which we tackle are beyond the area of
our professional interest and we do not claim any thorough
understanding of them still, the questions formulated in the language
of our trade seem to be extremely interesting for its own,
mathematical meaning.  In fact, we believe that most of the important
things in mathematics are related to real phenomena of Nature. The
interpretation of this profound feature of Mathematics is left for the
reader and it will definitely depend on the reader's attitude towards
fundamental Creation vs. Evolution problem, cf.~\cite{shafarevich} and
\cite{reid1}.

Knowing our limitations as laymen in computational biology and
statistics we try to stay within borders of the branch of
mathematics which we believe we understand. That is why we take a
relatively simple model, redefine it in purely algebraic language and
examine it using methods of algebraic geometry. The result exceeds our
original expectations, we find the object appearing in this process
very interesting for its own, pure geometric aspects, with properties which
we have not expected originally.

Our original task was computing Hilbert-Ehrhard polynomials for
varieties arising as geometrical models of binary symmetric 3-valent
phylogenetic trees. The question is consistent with the attitude of
computational algebraic geometry and algebraic statistics where the
point is to compute and understand the ideal of the variety in
question in the ambient projective space. Then the Hilbert-Ehrhard
polynomial provides a fundamental invariant of such an ideal, the
dimensions of homogeneous parts of it. To our surprise the polynomial
does not depend on the shape of the tree but merely on its size, the
number of leaves or, equivalently the dimension of its geometric
model. The strive to understanding this phenomenon lead us to proving
one of the main results of the present paper,
\ref{models-deformation-equivalent}, which asserts that models of
trees with the same number of leaves are deformation equivalent, that
is they are in the same connected component of the Hilbert scheme of
the projective space in question (hence they have the same Hilbert
polynomial).

The fact that the geometric models of trees modelling some processes
--- the discreet objects --- live in a connected continuous family of
geometric objects probably deserves its explanation in terms of
algebraic statistic or even biology. For the algebraic geometry part
we have a natural question arising about irreducibility of the
component of the Hilbert scheme containing these models and (if the
irreducibility is confirmed) about varieties which arise as general
deformations (that is, over a general point of the component of the
Hilbert scheme in question).  The question about a general deformation
of the model is related to the other main result of the present paper,
\ref{Fano}, which is that these models are index 4 Fano varieties with
Gorenstein teminal singularities. Thus one would expect that their
general deformation is a {\em smooth} Fano variety of index 4,
c.f.~\cite{namikawa}.

The present paper is organized as follows. We deal with varieties
defined over complex numbers. In the first section we define
phylogenetic trees and their geometric models. We do it in pure
algebraic way and with many simplifications: we deal with unrooted
symmetric trees which are then assumed to be binary and eventually
3-valent. From the algebraic geometer point of view studying geometric
models in this case can be reduced to understanding special linear
subsystems of the Segre linear system on a product of $\P^1$'s,
\ref{complete-model=Segre}. Eventually, the question boils down to
studying fixed points of the Segre system with respect to an action of
a group of involutions, \ref{map-injective}. Since the action can be
diagonalized this brings us down to toric geometry.

In the second section we define a geometric model of a tree in terms
of toric geometry, via a polytope in the space of characters of a
complex torus, which we call a polytope model of the tree and to which
we subsequently associate a projective variety. The main results of
this part are \ref{polytope=>parametrization} and
\ref{equivalence-of-models} which assert that the models defined in
the first part are the same as these defined the toric way. In this
part we also prove results which are of the fundamental technical
importance: this is a fiber product formula for polytopes of trees,
\ref{polytope-of-graft-is-fiberproduct}, and its counterpart for
varieties, a quotient formula \ref{thm-graft-as-quotient}. The latter
asserts that the geometric model of a tree obtained by gluing two
smaller trees is a Mumford's GIT (Geometric Invariant Theory) quotient
of the product of their respective models.

The third section of the present paper contains its main
results. After a brief discussion of equations defining a geometric
model of a tree, with special consideration to a tree with two inner
nodes and four leaves, we examine fans of geometric models and
resolution of their singularities. We prove that geometric models of
3-valent binary symmetric trees are index 4 Fano varieties with
Gorenstein terminal singularities which admit small resolution,
\ref{Fano}. Next we consider deformations of models of trees. The
approach is, roughly, as follows: we know how to deform equations of a
small tree with four leaves and one inner edge, the result of the
deformation is another tree with the inner edge ``flopped'':
$$
\begin{array}{cccccccc}
\begin{xy}<7pt,0pt>:
(-1,0)*={}="0" ; (-2,2)*={1\phantom{..}} **@{-},
"0" ; (-2,-2)*={2\phantom{.}} **@{-},
"0" ; (1,0)*={}="1" **@{-},
"1" ; (2,2)*={\phantom{.}3} **@{-},
"1" ; (2,-2)*={\phantom{.}4} **@{-},
\end{xy}
&&\longleftrightarrow&&
\begin{xy}<7pt,0pt>:
(0,1)*={}="0" ; (-2,2)*={1\phantom{..}} **@{-},
"0" ; (0,-1)*={}="1" **@{-},
"1" ; (-2,-2)*={2\phantom{.}} **@{-},
"0" ; (2,2)*={\phantom{.}3} **@{-},
"1" ; (2,-2)*={\phantom{.}4} **@{-},
\end{xy}
\end{array}
$$ Applying the GIT quotient formula, \ref{thm-graft-as-quotient}, we
are able to use this elementary deformation associated to four leaves
trees to get a similar deformation for {\em every} inner edge of any
tree, \ref{quotient-of-product}. This implies the result about
deforming one geometric model to another,
\ref{models-deformation-equivalent}.

In the last part of section 3, we discuss Hilbert-Ehrhard polynomial
of models (both polytopes and varieties) of trees. We define a
relative version of the polynomial and then a product of such
polynomials which is related to gluing respective trees. The
elementary deformation procedure implies associativity of the product
which not only implies the invariance of the Hilbert-Ehrhard
polynomial for trees with the same number of leaves but also provides
a simple formula for computing it, \ref{relative-ehrhard-equal}.

The appendix contains some computations. Firstly we prove that the
polytopes of the 3-valent trees are normal which is needed to ensure
the proper definition of their geometrical models. Next, using the
\cite{polymake} software we verify a simple (yet 9-dimensional)
example to check that the polytope models of different trees in this
case are different. The question if the polytope (or geometric) models
of non-isomorphic trees are non-isomorphic is open,
c.f.~\cite{allman-rhodes2}. Finally, we make numerical experiments
(using \cite{maxima} and \cite{gnuplot}) to look at the behaviour of
the relative volume distribution which measures the (normalized) volume
of the model with respect to a fixed leaf of a tree.

The paper uses consistently the language of algebraic geometry,
including toric geometry. We ignore, or barely mention, relations to
algebraic statistic and biology, suggesting the reader to look into
\cite{ps} (or into \cite{erss} for a concise version of exposition),
to get an idea about the background of the problems that we deal
with. It was our primary intention to make the present paper
self-contained so that it can be read as it is by an algebraic
geometer with no knowledge of its possible applications outside
algebraic geometry. On the other hand, a reader who is not familiar
with algebraic geometry but is interested in acquiring ideas which are
important in our approach (regarding quotients and deformations) is
advised to look into \cite{reid2} and \cite{altman} for a short
exposition to these matters.

We would like to thank Jaros{\l}aw Buczy{\'n}ski for his remarks and
Piotr Zwiernik for bringing this subject to our attention.

\subsection{Notation}

\renewcommand{\labelitemi}{$\circ$}

\begin{itemize}

\item $|\cA|$ denotes cardinality of a finite set $\cA$.

\item A lattice is a finitely generated free abelian group.

\item Depending on the context a subscript denotes the extension of
the basic ring or a fiber of a morphism, e.g.~$M_\R=M\otimes_\Z\R$.

\item Given a finite dimensional vector space (or a lattice) $V$ with
a basis $\{ v_1,\dots, v_n\}$, by $\{v_1^*,\dots,v_n^*\}$ we will
denote the dual base of $V^*$, that is $v^*_i(v_i)=1$ and
$v^*_i(v_j)=0$ if $i\ne j$.

\end{itemize}


\section{Preliminaries: phylogenetic trees.}


{\bf Summary:} (for algebraic geometers) phylogenetic trees are a
clever way of describing linear subsystems of Segre system on the
product of projective spaces. In case of binary symmetric trees the
question is to find subsystems of sections of Segre system on a
product of $\P^1$'s invariant with respect to some $\Z_2^{|\N|}$
action.

\subsection{Trees and linear algebra}

\begin{notation}
A tree $\T$ is a simply connected graph (1-dimensional CW complex)
with a set of edges $\E=\E(\T)$ and vertices $\V=\V(\T)$ and the
(unordered) boundary map $\d: \E \to \V^{\wedge 2}$, where $\V^{\wedge
2}$ denotes the set of unordered pairs of distinct elements in
$\V$. The number $|\E|\geq 1$ is, by definition, the number of edges
of $\T$, then number of vertices $|\V|$ is $|\E|+1$.  We write
$\d(e)=\{\d_1(e),\d_2(e)\}$ and say $v$ is a vertex of $e$, or $e$
contains $v$ if $v\in\{\d_1(e), \d_2(e)\}$, we simply write $v\in e$.
The valency of a vertex $v$ is the number of edges which contain $v$
(the valency is positive since $\T$ is connected and we assume it has
at least one edge). A vertex $v$ is called a leaf if its valency is 1,
otherwise it is called an inner vertex or a node.  If the valency of
each inner node is $m$ then the tree will be called $m$-valent. The
set of leaves and nodes will be denoted $\cL$ and $\N$, respectively,
$\V=\cL\cup\N$.  An edge which contains a leaf is called a petiole, an
edge which is not a petiole is called an inner edge (or branch), and
the set of inner edges will be denoted by $\E^o$.
\end{notation}

\begin{ex}
An {\em caterpillar} of length $n$ is a 3-valent tree with $n$ inner
edges and $n+1$ inner nodes whose defoliation (i.e.~after removing all
leaves and petioles) is just a string of edges. That is, there are
exactly two inner nodes to which of them there are attached two
petioles (we call them heads or tails), any other inner node has
exactly one petiole (called a leg) attached.

$$
\begin{xy}<15pt,0pt>:
(0,-0.1)*={}="0" ; (-\halfroottwo,\halfroottwo)*={} **@{-},
"0" ; (-\halfroottwo,-\halfroottwo)*={} **@{-},
"0" ; (1,0.1)*={}="1" **@{-},
"1" ; (0.9,-1)*={}     **@{-},
"1" ; (2,0.1)*={}="2" **@{-},
"2" ; (2.1,-1)*={}     **@{-},
"2" ; (4,0)*={}="4" **@{.},
"4" ; (3.9,-1)*={} **@{-},
"4" ; (5,0.1)*={}="5"  **@{-},
"5" ; (\fiveandhalfroottwo,\halfroottwo)*={} **@{-},
"5" ; (\fiveandhalfroottwo,-\halfroottwo)*={} **@{-},
\end{xy}
$$

\end{ex}

\begin{notation}
Let $W$ be a (complex, finite dimensional) vector space with a
distinguished basis, sometimes called letters: $\{ \a_0, \a_1,
\a_2\dots\}$.    We consider the
map $\suma: W \to \C$, such that $\suma(\a_i)=1$ for every $i$, that
is $\suma=\sum\a^*_i$.

Let $\hW$ be a subspace of the second tensor product $W\otimes W$.  An
element $\sum_{i,j}a_{ij}(\a_i\otimes\a_j)$ of $\hW$ can be represented
as a matrix $\left(a_{ij}\right)$. Through the present paper we will
assume that these matrices are symmetric, that is $\hW$ is contained in
$S^2(W)$.

Given a tree $\T$ and a vector space $W$, and a subspace $\hW\subset
W\otimes W$ we associate to any vertex $v$ of $\V(\T)$ a copy of $W$
denoted by $W_v$ and for any edge $e\in\E(\T)$ we associate a copy of
$\hW$ understood as the subspace in the tensor product $\hW^e\subset
W_{\d_1(e)}\otimes W_{\d_2(e)}$. Note that although the pair
$\{\d_1(e),\d_2(e)\}$ is unordered, this definition makes sense since
$\hW$ consists of symmetric tensors. Elements of $\hW^e$ will be
written as (symmetric) matrices $(a^e_{\a_i,\a_j})$.
\end{notation}

\begin{defin}\label{def-phylogenic-tree}
The triple $(\T,W,\hW)$ together with the above association is called
a (symmetric, unrooted) phylogenetic tree.
\end{defin}

\begin{const}\label{parameterization}
Let us consider a linear map of tensor products
$$
\widehat\Psi:\hW^\E=\bigotimes_{e\in\E}\  \hW^e \lra 
W_\V = \bigotimes_{v \in \V} W_v
$$
defined by setting its dual as follows
$$
\widehat\Psi^*(\otimes_{v\in\V} \ \ \a_{v}^*)= 
\otimes_{e\in\E}\ (\a_{\d_1(e)}\otimes\a_{\d_2(e)})^*_{|\hW^e}
$$ where $\alpha_v$ stands for an element of the chosen basis
$\{\alpha_i\}$ of the space $W_v$.  The complete affine geometric model of
the phylogenetic tree $(\T,W,\hW)$ is the image of the associated
multi-linear map
$$
\widetilde\Psi: \product_{e\in\E}\hW^e \lra
W_\V=\bigotimes_{v \in \V} W_v
$$
The induced rational map of projective varieties will be denoted by $\Psi$: 
$$
\Psi: \product_{e\in\E} \  \P(\hW^e) \ \ratmap \P(W_\V)=
\P(\bigotimes_{v \in \V} W_v)
$$
and the closure of the image of $\Psi$ is called the complete projective
geometric model, or just the complete model of $(\T,W,\hW)$. The maps
$\widetilde\Psi$ and $\Psi$ are called the parameterization of the
respective model.

Given a set of vertices of the tree we can ``hide'' them by applying
the map $\suma=\sum_i\a^*_i$ to their tensor factors. In what follows
will hide inner nodes and project to leaves. That is, we consider the
map
$$
\begin{array}{c}
\Pi_{\L}:
W_\V= \bigotimes_{v\in\V} W_v \to W_\L=\bigotimes_{v\in\L} W_v\\ \\
\Pi_{\L}=(\otimes_{v\in\L}\ id_{W_v}) \otimes(\otimes_{v\in\N}\
\suma_{W_v})
\end{array}
$$ 
\end{const}

\begin{defin}\label{def-model}
The affine geometrical model of a phylogenetic tree $(\T,W,\hW)$ is an
affine subvariety of $W_\L=\bigotimes_{v\in\L} W_v$ which is the image
of the composition $\Phi=\Pi_{\L}\circ\Psi$. Respectively, the
projective geometrical model, or just a model, denoted by $X(\T)$ is
the underlying projective variety in $\P(W_\L)$. For $X=X(\T)$ by
$\cO_X(1)$ we will denote the the hyperplane section bundle coming
from the embedding in the projective space $\P(W_\L)$.
\end{defin}

Note that $X(\T)$ is the closure of the image of the respective rational map 
$$
\product_{e\in\E}\P(\hW^e)\ \ratmap \P\left(\bigotimes_{v\in\L} W_v\right)
$$ 
which is defined by a special linear subsystem in the Segre linear system
$|\bigotimes_{e\in\E} p^*_{\P(\hW^e)}\cO_{\P(\hW^e)}(1)|$, where
$p^*_{\P(\hW^e)}$ is the projection from the product to the respective
component. We will call this map a {\em rational parametrization} of
the model.

The above definition of parametrization is an unrooted and algebracized
version of what is commonly considered in the literature, see 
e.g.~\cite{allman-rhodes1}, \cite{ss} or \cite{cgs}.

\subsection{Binary symmetric trees.}\label{bin-symm-trees}


Depending on the choice of $\hW \subset W\otimes W$ we get different
phylogenetic trees and their models.  A natural assumption is that in
the matrix representation the elements of $\hW$ the sum of the numbers
in each row and each column is the same (in applications, these
numbers would stand for the probability distribution so their sum
should be equal to 1).  If $W$ is of dimension 2 this is equivalent to
saying that the respective matrix is of the form
$$\left[\begin{array}{cc} a&b\\ b&a \end{array}\right]$$ for some
$a$ and $b$ in $\C$.

From now on we will consider binary symmetric phylogenetic trees, that
is, we assume that dimension of $W$ is 2 and $\hW$ consists of
matrices (tensors) satisfying the above symmetric condition. The
elements of the distinguished basis of $W$ will be denoted $\a$ and
$\b$. Note that $\hW$ has dimension 2 as well. We will call them
binary symmetric trees or just trees when the context is obvious. Our
task is to understand geometric models of these trees.

\begin{ex}\label{3startree-map}
Let $\T$ be a tree which has one inner node $v_0$, three leaves $v_1,\
v_2,\ v_3$ whose petioles we denote, respectively, by $e_1,\ e_2,\
e_3$. We denote the basis of $W_{v_i}$ by $\alpha_i, \ \beta_i$, while
$\widehat{W}^{e_i}$ consists of matrices $\left[\begin{array}{cc}
a_i&b_i\\ b_i&a_i \end{array}\right]$. Then the parameterization map
$$\widetilde\Psi: \hW^{e_1}\times\hW^{e_1}\times\hW^{e_3} \lra
W_{v_1}\otimes W_{v_2}\otimes W_{v_3}$$
is as follows:
$$\begin{array}{l}
\widetilde\Psi(a_1,b_1,a_2,b_2,a_3,b_3)=\\
(a_1a_2a_3+b_1b_2b_3)\cdot
(\alpha_1\otimes\alpha_2\otimes\alpha_3+\beta_1\otimes\beta_2\otimes\beta_3)\  
+\\
(b_1a_2a_3+a_1b_2b_3)\cdot
(\beta_1\otimes\alpha_2\otimes\alpha_3+\alpha_1\otimes\beta_2\otimes\beta_3)\  
+\\
(a_1b_2a_3+b_1a_2b_3)\cdot
(\alpha_1\otimes\beta_2\otimes\alpha_3+\beta_1\otimes\alpha_2\otimes\beta_3)\  
+\\
(a_1a_2b_3+b_1b_2a_3)\cdot
(\alpha_1\otimes\alpha_2\otimes\beta_3+\beta_1\otimes\beta_2\otimes\alpha_3)\\
\end{array}
$$
\end{ex}


\begin{notation}
Let $\rho: W \to W$ be a linear involution $\rho(\a)=\b$,
$\rho(\b)=\a$, the map $\rho$ is reflection with respect to the linear
space $W^\rho$ spanned by $\a+\b$.
We note that on $\hW$ the right and left
action of $\rho$ coincide, i.e. 
$(\rho\otimes id_W)_{|\hW}=(id_W\otimes\rho)_{|\hW}$,
and the resulting involution will be denoted by $\widehat\rho$, note that
$$\widehat\rho\left(\left[\begin{array}{cc}a&b\\ b&a\end{array}\right]\right) =
\left[\begin{array}{cc}b&a\\ a&b\end{array}\right]$$  
In particular,
$\rho\otimes\rho$ is identity on $\hW$.

Given a binary symmetric tree $(\T,W,\hW)$ we define respective
involutions: $$\rho_\V=\otimes_{v\in\V} \ {\rho_v}:
\bigotimes_{v\in\V} W_v \to \bigotimes_{v\in\V} W_v$$
$$\rho_\L=\otimes_{v\in\L} \ \rho_v: \bigotimes_{v\in\L} W_v \to
\bigotimes_{v\in\L} W_v$$ Let $W_\V^{\rho}=(\bigotimes_{v\in\V}
W_v)^{\rho_\V}$ and $W_\L^{\rho}=(\bigotimes_{v\in\L} W_v)^{\rho_\L}$
be their fixed points, that is the maximal subspace on which $\rho_\V$
and, respectively, $\rho_\L$ acts trivially.
\end{notation}

\begin{lemma}\label{complete-model=Segre}
The image of $\Psi$ is contained in $\P(W_\V^{\rho})$ and
the induced map $$\Psi:\product_{e\in\E}\P(\hW^e) \to
\P(W_\V^{\rho})$$ is Segre embedding.
\end{lemma}

\begin{pf} We want to prove that $\widehat\Psi$ maps 
$\hW^\E$ isomorphically to the space $W_\V^{\rho}$.
First let us note that
$$\begin{array}{c} \widehat\Psi^*(\rho^*_\V(\otimes_{v\in\V}\ \a_v^*))=
\widehat\Psi^*(\otimes_{v\in\V}\ \rho^*_v(\a_v^*))= \\ \otimes_{e\in\E}
(\rho^*_{\d_1(e)}(\a_{\d_1(e)}^*)\otimes\ \rho^*_{\d_2(e)}(\a_{\d_2(e)}^*))=
\otimes_{e\in\E} (\a_{\d_1(e)}^*\otimes\a_{\d_2(e)}^*)=
\widehat\Psi^*(\otimes_{v\in\V}\ \a_v^*)
\end{array}
$$ (where, again, as in \ref{parameterization} $\alpha_v$ denotes
either $\alpha$ or $\beta$ in the space $W_v$) so that
$\widehat\Psi^*\circ\rho^*_\V=\widehat\Psi^*$ which implies
$\rho_\V\circ\widehat\Psi=\widehat\Psi$ , hence
$\im(\widehat\Psi)\subset W_\V^\rho$.

Next, let us note that $\dim W_\V^\rho=2^{|\V|-1}$ so that it is equal
to $\dim \hW^\E$ because $|\E|=|\V|-1$. The proof (e.g.~by induction
with respect to $|\V|$) is instantaneous if one observes that the
basis of $W_\V$ can be made of tensor products of $(+1)$ and $(-1)$
eigenvectors of each $\rho_v$ and thus $W_\V$ splits into the sum of
$(+1)$ and $(-1)$ eigenspaces of $\rho_\V$, each of the same
dimension.

Now, to conclude the proof we have to show that $\widehat\Psi$ is
injective which is equivalent to $\widehat\Psi^*$ being
surjective. Note that $\hW^*$ is spanned by two forms:
$$\begin{array}{cc}
\gamma_0\left(\left[\begin{array}{cc}a&b\\ b&a\end{array}\right]\right)=a,
&
\gamma_1\left(\left[\begin{array}{cc}a&b\\ b&a\end{array}\right]\right)=b
\end{array}$$ 
and $\gamma_1=\gamma_0\circ\widehat{\rho}$.

Now given an element $\otimes_{e\in\E}\ \gamma_{i(e)}\in \bigotimes_{e\in\E}
\hW^e$ we define inductively a sequence $\a_v^*$, indexed by vertices
of $\T$ such that $\widehat\Psi(\otimes_{v\in\V}\ \a_v^*)=
\otimes_{e\in\E}\ \gamma_{i(e)}$. We choose a vertex $v_0$ and
set $\a^*_{v_0}$ to be either $\a^*$ or $\beta^*$. Now suppose that
$\a^*_{v}$ are defined for $v$ in a subtree $\T'$ of $\T$. Suppose
that $v'$ is not in $\T'$ but is joined to a vertex $v''$ in $\T'$ by
and edge $e'$.  Then we set $\a^*_{v'}=\a^*_{v''}$ if
$\gamma_{i(e')}=\gamma_0$ or $\a^*_{v'}=\rho(\a^*_{v''})$ if
$\gamma_{i(e')}=\gamma_1$.

\end{pf}

\begin{notation}
Let us choose a node $v\in\N$.  We consider an involution
$\widehat{\rho}_v^\E$ on the space
$\hW^\E=\bigotimes_{e\in\E}\hW^e$. First, for any $e\in\E$ we set
$\widehat\rho_v^e=id_{\hW^e}$ if $v$ is not a vertex of $e$ and
$\widehat\rho_v^e=\widehat\rho$ if $v$ is a vertex of $e$. Next, we
define $\widehat\rho^\E_v=\bigotimes_{e\in\E}\widehat\rho^e_v$.  Let
$G_\N$ be the group of automorphisms of $\hW^\E$ generated by
involutions $\rho_v^\E$, for $v\in\N$.
\end{notation}

The following observation about a convenient choice of coordinates is
sometimes referred to as a Fourier transform, see e.g.~\cite{ss}.

\begin{lemma}\label{diagonalization}
$G_\N\iso\Z_2^{|\N|}$, the action of $G_\N$ restricts to
$\product_{e\in\E}\hW^e$ and it is is diagonalizable. 
\end{lemma}

\begin{pf}
This follows immediately from the definition of the action of
$\widehat\rho_v^e$ on $\hW^e$. Namely, in basis of $\hW^e$ consisting of $a+b$
and $a-b$ the action of $\widehat\rho_v^e$ is diagonal. 
\end{pf}

\begin{lemma}\label{map-injective}
The map $\widehat\Psi^*\circ\Pi^*_\cL$ maps
$(W_\L^\rho)^*$ invectively into the space $((\hW^\E)^*)^{G_\N}$ 
on which $G_\N$ acts trivially.
\end{lemma}

\begin{pf}
The proof is similar to that of \ref{complete-model=Segre}. Firstly,
the map $\hat\Psi^*\circ\Pi^*_\cL$ is injective on $(W_\cL^\rho)^*$
because of \ref{complete-model=Segre} and injectivity of $\Pi^*_\cL$.
Next, we note that the action of $G_\N$ is trivial on its image.
Indeed, we define
$\rho^\V_v=\bigotimes_{w\in\V}\rho^w_v$, where, for any
$w\in\V$ we set $\widehat\rho_v^w=id_{W_w}$ if $w\ne v$ and
$\rho_v^w=\rho$ if $w=v$. Then for $v\in\N$ we have
$$\begin{array}{cccccc}
(\widehat\rho^\E_v)^*\circ\widehat\Psi^*=\widehat\Psi^*\circ(\rho^\V_v)^*
&{\rm and}&
\Pi_\cL\circ\rho^\V_v=\Pi_\cL
\end{array}
$$ where the first equality follows directly from the definition of
the map $\widehat\Psi$, \ref{parameterization}. This implies
$(\widehat\rho^\E_v)^*\circ\widehat\Psi^*\circ\Pi^*_\cL=
\widehat\Psi^*\circ\Pi^*_\cL$
which is what we want.
\end{pf}

We will prove that, in fact, $\widehat\Psi^*\circ\Pi^*_\cL$ is an
isomorphism, \ref{polytope=>parametrization}, so that the geometric
model of the tree is defined by $G^\N$ invariant sections of the Segre
linear system. Because of \ref{diagonalization} $G^\N$ can be treated
as a subgroup of a complex torus and thus we can use toric geometry.


\section{Toric geometry.}

{\bf Summary:} We study invariants of an action of $\Z_2^\N$ on
$(\P^1)^{\times |\E|}$ and a related polytope in the cube
$[0,1]^{|\E|}$ which we call a polytope model of the tree.  The
polytope models are used to define geometric models in terms of toric
geometry. These polytopes turn out to be fiber products of elementary
ones. This leads to interpreting the geometrical model of a tree as a
quotient of products.


\subsection{Lattice of a tree and the action of the torus}\label{sect-torus}


Given a tree $\T$ we encode it in terms of dual
lattices.

\begin{defin}\label{def-lattice-of-tree}
Let $\T$ be a tree with the set of vertices $\V$ and the set of edges
$\E$. We set $M=M(\T)=\bigoplus_{e\in\E} \Z\cdot e$ to be a lattice,
or free abelian group, spanned on the set $\E$. Let
$N=N(\T)=Hom(M,\Z)$ be the dual lattice.  We represent elements of
$\V$ as elements of $N$. Namely, for $v\in\V$ we set $v(e)=1$ if $e$
contains the vertex $v$ and $v(e)=0$ otherwise. The pair $(M,N)$
together with the choice of the basis $\E$ of $M$ and set
$\V\subset N$ is called the lattice pair of the tree $\T$.
\end{defin}

From this point on we identify the edges and the vertices of $\T$ with
the respective elements in $M(\T)$ and $N(\T)$.  The elements of the
basis of $N$ dual to $\{e\in\E\}$ will be denoted by $e^*$.  Then for
any $v\in\V$ we have, by definition, $v=\sum_{e\ni v} e^*: N\ra
\Z$. In particular, $v$ is a leaf if and only if $v=e^*$ for some $e$,
which is a petiole for $v$.  

Let us recall that $|\V|=|\E|+1$ so the set of vertices has to be
linearly dependent in $N$. The set of vertices of $\T$ can be divided
into two disjoint classes, say $\V=\V^-\cup\V^+$, each class
consisting of vertices which can be reached one from another by
passing through an even number of edges.

\begin{lemma}\label{linear-relation-in-tree}
The equality $\sum_{v\in\V^-}v=\sum_{v\in\V^+}v$ is, up to
multiplication by a constant, the only linear relation in $N$ between
vectors $v$ from $\V$. In particular, any proper subset of $\V$
consists of linearly independent vectors in $N$.
\end{lemma}

\begin{pf} Suppose that $\sum a_v\cdot v=0$, for some $a_v\in k$. 
For any $e\in\E$ we have 
$$(\sum a_v\cdot v)(e)= a_{\partial_1(e)}+a_{\partial_2(e)}$$ and
therefore $a_{\partial_1(e)}=-a_{\partial_2(e)}$. Thus we get the
desired relation.
\end{pf}
The operations on trees can be translated to lattices, here is an
example.
\begin{const}
Let $v_0$ be a
2-valent inner node of $\T$ which belongs to exactly two edges $e_1$
and $e_2$. Let $\T_{v_0}$ be a tree obtained by removing the node
$v_0$ from $\T$ and replacing the edges $e_1$ and $e_2$ by a single
edge $e_0$.
Let $(M,N)$ be the lattice pair of $\T$. We set $M_{v_0}\subset M$ to be
the kernel of $e^*_2-e^*_1$ and $N_{v_0}=N/\Z\cdot(e^*_2-e^*_1)$, clearly
$M_{v_0}$ and $N_{v_0}$ are dual. We define $e_0'=e_1+e_2$. Note that
$\E_{v_0}=\E\setminus\{e_1,e_2\}\cup\{e_0'\}$ is a basis of $M_{v_0}$.  For
$v\in\V\setminus\{v_0\}$ by $v'$ we denote the image of a vertex $v$
under the projection $N\ra N_{v_0}$ and set of all $v'$ we denote by
$\V_{v_0}$. One can verify easily the following.
\end{const}

\begin{lemma}\label{lattice-removing2valentnode}
The above defined pair $(M_{v_0},N_{v_0})$ together with the above
choice of $\E_{v_0}$ and $\V_{v_0}$ is the lattice pair of the tree
$\T_{v_0}$ obtained from $\T$ by removing the 2-valent inner node $v_0$.
\end{lemma}

Now we set up the toric environment.

\begin{const}
We deal with a binary symmetric tree $(\T,W,\hW)$.  Because of
\ref{diagonalization} for any edge $e\in\E$ there exists an
inhomogeneous coordinate $z_e$ on $\P^1_e=\P(\hW^e)$ such that for
$v\in e$ the action of $\rho_v^\E$ is as follows
$\rho_v^\E(z_e)=-z_e$. 

Let $T\iso\product_{e\in\E}\C^*$ be a torus with
coordinates $\{z_e\in\C^*: e\in\E\}$ and with the natural action
$$T\times \product_{e\in\E}\P^1_e\lra \product_{e\in\E}\P^1_e$$ which
is the multiplication of $z_e$'s coordinate-wise.  We consider an
injective map $\iota: G_{\N}\to T\iso\product_{e\in\E}\C^*$ which is
defined as follows. For any $\rho_v^\E\in G_{\N}$ we take
$\iota(\rho_v^\E)\in T$ such that $z_e(\iota(\rho_v^\E))=-1$ if $v\in
e$ and $z_e(\iota(\rho_v^\E))=1$ if $v\not\in e$. Then $\iota$ extends
to a homomorphism of groups $\iota: G_{\N}\to
T\iso\product_{e\in\E}\C^*$ and the action of $G_{\N}$ on
$\product_{e\in\E}\P^1_e$ factors through $\iota$.

We explain this situation using the lattices of the tree $\T$ and
toric geometry formalism. Our notation is consistent with this of
standard toric geometry textbooks, e.g.~\cite{oda} or \cite{fulton}.
We take the torus $T=T_N=N\otimes_\Z\C^*$ with coordinates
$z_e=\chi^e$, where $e\in\E\subset M$ is the distinguished
basis. Recall that the elements of $M$ can be identified with monomials
in coordinates $z_e$, that is, each $u\in M$ such that $u=\sum a_ie_i$
represents a monomial $\chi^u=\product z_{e_i}^{a_i}$. For $w\in N$
and $t\in\C^*$ the $z_e$-th coordinate of the respective point
on $T_N=N\otimes_\Z\C^*$ is as follows
$z_e(w\otimes t)=t^{w(e)}$.
Moreover, recall that every element $w$ of $N$ can be identified with
algebraic 1-parameter subgroups $\lambda_w$ of $T_N$. That is, for
$w\in N$ and $t\in\C^*$ we set
$\lambda_w(t)(z_e)=t^{w(e)}\cdot z_e$.  In short,
$N=\Hom_{alg}(\C^*,T_N)$ and $M=\Hom_{alg}(T_N,\C^*)$, \cite[Sect 2.3]{fulton}

The complexified lattice $N_\C=N\otimes_\Z\C$ can interpreted as the
tangent space to the unit element in the torus $T_N$ and we have the
natural exponential map $N_\C\ra T_N$. In particular, $z_e(\exp(2\pi i
(w)))=\exp(2\pi i(w(e))$. The image of the real vector space
$N_\R\subset N_\C$ under this exponential map is the maximal compact
real subgroup $\prod S^1$ of $T_N$. Using the exponential map we can
relate the vertices $v\in\N$ viewed as elements of the lattice $N$ to
their respective automorphisms $\rho^\E_v\in G^\N$.  The following
lemma is not used directly in our arguments so we skip its proof.
\end{const}

\begin{lemma}\label{diagonalization2}
For every $v\in\V$ we have $exp(2\pi i(v/2))=\iota(\rho_v^\E)$.  If
$\widehat N$ is the lattice spanned in $N_\R$ by $N$ and $\N/2=\{v/2:
v\in\N\}$ then the inclusion $N\hookrightarrow\widehat N$ yields
an exact sequence of groups
$$0\lra \iota(G^\N)\lra T_N=N\otimes_\Z\C^*\lra 
\widehat N\otimes_\Z \C^*\lra 0$$
\end{lemma}

For our purposes we need the following lemma which provides a clear
description of functions on the torus $T_N$ which are invariant with
respect to the action of $G_{\N}$.

\begin{lemma}\label{invariant=even}
A monomial function $\chi^u$ on $T_N$ with $u\in M$ is invariant with respect
to the action of $\rho_v^\E$ if and only if $v(u)\in 2\Z$.
\end{lemma}

\begin{pf}
First, we note that, by definition,
$\rho_v^\E(\chi^e)=\rho_v^\E(z_e)=(-1)^{v(e)}\cdot z_e=
(-1)^{v(e)}\chi^e$. Next, we write the exponent of the monomial $\chi^u$ in
terms of the distinguished coordinates:
$u=\sum_{e\in\E}e^*(u)e$. Then, since $\rho_v^\E$ is a homomorphism, we get
$$\rho_v^\E(\chi^u)=(-1)^{\sum_e e^*(u)\cdot v(e)}\cdot\chi^u =
(-1)^{v(u)}\cdot \chi^u$$ which concludes the proof.
\end{pf}

\begin{defin}\label{def-normalized-lattice}
Given a tree $\T$ with the lattice pair $(M,N)=(M(\T),N(\T))$ we define
its normalized lattice pair $(\hM,\hN)=(\hM(\T),\hN(\T))$ as follows:
$\hM=\{u\in M:\forall v\in\V\ v(u)\in 2\Z\}$ and $\hN$ is a dual of $\hM$
which contains $N$ and the set $\N/2=\{v/2: v\in \N\}$.
\end{defin}

In view of \ref{invariant=even} the lattice $\widehat{M}$ contains
monomials which are $G^\N$ invariant.


\subsection{Polytope model of a tree}\label{sect-polytope-model}


The complete Segre linear system on $\product_{e\in\E}\P^1_e$ is
spanned on monomials $\product_{e\in\E} z_e^{\epsilon_e}$ where
$\epsilon_e\in\{0,\ 1\}$. Equivalently, once the big torus action
$T_N$ on $\product\P^1_e$ is chosen, the complete Segre system is
represented by vertices of the unit cube $\boxdot_M=\{u\in M_\R:
\forall_i\ 0\leq e^*_i(u)\leq 1\}$ in the space of characters $M_\R$,
or by zero-one sequences indexed by $\E$. 

Because of \ref{complete-model=Segre} and \ref{map-injective} we are
interested in subsystems of the Segre linear systems or,
equivalently, subsets of of vertices of $\boxdot_M$.  If $\Delta$ is a
polytope in $M_\R$ whose vertices are contained in the set of vertices
of $\boxdot_M$ then we call it a subcube.

\begin{defin}\label{polytope-model}
Given a binary tree $\T$ with its lattice pair $(M,N)$ its polytope
model $\Delta(\T)$ is a polytope in the lattice $M$ which is the
convex hull of $\{u=\sum a_ie_i\in M: a_i=0,\, 1 \, {\rm and}\,
v(u)\in2\Z\ \ {\rm for\, every}\, v\in\N\}$. 
\end{defin}

We note that the vertices of $\Delta$ are precisely these among
vertices of $\boxdot_M$ which are in the sublattice $\hM\subset M$ and
because of \ref{invariant=even} they are exactly these monomials in
the complete Segre system which are invariant with respect to the
action of $G_\N$.

Since the cube $\boxdot_M$ is the fundamental domain in dividing $M$
modulo 2 we can interpret the elements of the complete Segre system as
points in the linear space $M_{\Z_2}=M\otimes_\Z\Z_2$.

\begin{lemma}\label{polytope-model-Z_2}
If the vertices of the cube $\boxdot_M$ are identified with the
points in linear space $M\otimes\Z_2$ then vertices of $\Delta(\T)$
form the linear subspace $\N^\perp\subset M\otimes\Z_2$ of zeros of
forms $v\in N\otimes\Z_2$, where $v\in\N(\T)$.
\end{lemma}

\begin{pf}
This is a restatement of \ref{invariant=even}.
\end{pf}

\begin{cor}\label{number-of-vertices}
The polytope $\Delta(\T)$ has $2^{|\L|-1}$ vertices. 
\end{cor}

\begin{pf}
We use \ref{polytope-model-Z_2}: by \ref{linear-relation-in-tree} the
elements $v$'s are linearly independent in $N\otimes\Z_2$ so dimension
of the space of their zeroes in $M\otimes\Z_2$ is
$|\E|-|\N|=|\L|-1$.
\end{pf}

Using the above information we can conclude identifying the linear
subsystem in the Segre system which defines the projective model of a
binary symmetric tree.

\begin{thm}\label{polytope=>parametrization}
In the situation of section \ref{bin-symm-trees} the map
$\widehat\Psi^*\circ\Pi^*_\L$ maps $(W_\L^\rho)^*$ isomorphically to
$((\hW^\E)^*)^{G_\N}$. In particular, in terms of the toric
coordinates on $\product_{e\in\E}\P(\hW^e)$ introduced in section
\ref{sect-torus}, the rational parametrization map
$$ 
\product_{e\in\E}\P(\hW^e)\ \ratmap\P\left(W^\rho_\cL\right)\subset
\P\left(\bigotimes_{v\in\L} W_v\right)
$$ 
is defined by elements of the Segre linear system on
$\product_{e\in\E}\P(\hW^e)$ which are associated to vertices of
$\Delta(\T)$.
\end{thm}

\begin{pf}
By the construction the vertices of $\Delta(\T)$ are these monomial in
the Segre system which are invariant with respect to the action of
$G^\N$. In other words they form a basis for $((\hW^\E)^*)^{G_\N}$.
In \ref{map-injective} have proved that the parametrization map
injects $(W_\L^\rho)^*$ into the space $((\hW^\E)^*)^{G_\N}$ and now
by \ref{number-of-vertices} they are of the same dimension so
this is an isomorphism.  
\end{pf}

Thus we have determined that studying projective geometric models of
binary symmetric trees is essentially equivalent to understanding
their polytopes. We start with the simplest, in fact trivial, example.

\begin{ex}
Let $\T$ be a tree consisting of two leaves, two petioles $e_1$ and
$e_2$, and one inner node $v_0$. Then $\Delta(\T)$ is spanned on $0$
and $e_1+e_2$.
\end{ex}

More generally we have the following result which extends
\ref{lattice-removing2valentnode}.

\begin{lemma}\label{removing2valentnode}
Suppose that $\T$ is a tree with a 2-valent node $v_0$.  Let
$\T_{v_0}$ be a tree obtained from $\T$ by removing $v_0$, as in the
situation of lemma \ref{lattice-removing2valentnode}.  Then, under the
natural inclusion $M(\T_{v_0})\subset M(\T)$ we have
$\Delta(\T_{v_0})=\Delta(\T)$ and $\hM(\T_{v_0})= \hM(\T)\cap M(\T_{v_0})_\R$
\end{lemma}

\begin{pf}
We use the notation of \ref{lattice-removing2valentnode}, in
particular $e_1$ and $e_2$ denote the edges containing $v_0$ and
$M(\T_{v_0})=M_{v_0}=\ker(e_2^*-e_1^*)$.  Note that the parity of the
node $v_0=e_1^*+e_2^*$ is equivalent to that of $e_2^*-e_1^*$, in
particular for $u\in M(\T_{v_0})=\ker(e_2^*-e_1^*)$ we have $v_0(u)\in
2\Z$.  Since $\N(\T)=\N(\T_{v_0}) \cup \{v_0\}$ the conditions defined
by $\N(\T)$ and $\N(\T_{v_0})$ on $M(\T_{v_0})_\R$ are the same.
Similarly, the conditions defining $\Delta(\T)$ and $\Delta(\T_{v_0})$
inside $M_{v_0}\otimes\R\cap \boxdot_M=\boxdot_{M_{v_0}}$ are the
same.
\end{pf}

By \ref{toric-via-polytope} we have that removing the 2-valent node
does not change the model of the tree. Thus, from now on we consider
trees with no 2-valent nodes.

\medskip
A star tree is a tree which has exactly one inner node, a star tree
with $d$ leaves will be denoted by $\T^{\!\!*d}$. 

\begin{lemma}\label{dim-star-trees-polytopes}
If $d\geq 3$ then vertices of $\Delta(\T^{\!\!*d})$ generate
$\hM(\T^{\!\!*d})$, in particular $\dim\Delta(\T^{\!\!*d})=d$.
\end{lemma}

\begin{pf}
If $\{e_i\}$ is the set of edges then $\Delta(\T^{\!\!*d})$ contains
sums $e_i+e_j$ for all possible pairs $i\ne j$ and for $d\geq 3$ they
span $\hM(\T^{\!\!*d})$ which is of index 2 in $M(\T^{\!\!*d})$.
\end{pf}

\begin{ex}\label{ex-tetrahedron}
The vertices of the polytope $\Delta(\T^{\!\!*3})$ with edges $e_0,\
e_1,\ e_2$ are as follows: $0$, $e_1+e_2$, $e_2+e_0$ and $e_0+e_1$ so
that $\Delta(\T^{\!\!*3})$ is a 3-dimensional tetrahedron. If
$\widehat M\subset M$ is the sublattice spanned by the vertices of
$\Delta(\T^{\!\!*3})$ then $M/\widehat M\iso\Z_2$. 

The inequalities defining $\Delta(\T^{\!\!*3})$ are as follows
$$
(-v/2)(\ \cdot \ )\geq -1, \ \ (v/2 - e_0^*)(\ \cdot \ )\geq 0, \ \ 
(v/2 - e_1^*)(\ \cdot \ )\geq 0, \ \ (v/2 - e_2^*)(\ \cdot \ )\geq 0.
$$ 
where, recall, $v=e_0^*+e_1^*+e_2^*$. 
\end{ex}

\begin{const}\label{const-fiber-product}
A {\em pointed tree} $(\T,\ell)$ is a pair consisting of a tree $\T$ and a
leaf $\ell\in\L(\T)$. Given two pointed trees $(\T_1,\ell_1)$ and
$(\T_2,\ell_2)$ we define their 
graft as follows:
$\T=\T_1\phantom{*}_{\ell_1}\!\!\!\vee_{\ell_2}\T_2$ is a tree
obtained by removing from each $\T_i$ the leaf $\ell_i$ and
identifying their respective petioles which becomes an inner edge
of the resulting tree $\T$. 

For example, a graft of two trees of type $\T^{\!\!*3}$ with
distinguished leaves denoted by $\circ$ is the following operation
$$
\begin{xy}<8pt,0pt>:
(0,0)*={}="1" ; (-0.5,\halfrootthree)*={} **@{-},
"1" ; (-0.5,-\halfrootthree)*={} **@{-},
"1" ; (0.8,0) **@{-},(1,0)*={\circ} 
\end{xy}\ 
\vee\ 
\begin{xy}<8pt,0pt>:
(0,0)*={}="2" ; (0.5,\halfrootthree)*={} **@{-},
"2" ; (0.5,-\halfrootthree)*={} **@{-},
"2" ; (-0.8,0) **@{-},(-1,0)*={\circ}
\end{xy}\ 
=\ 
\begin{xy}<8pt,0pt>:
(0,0)*={}="0" ; (-0.5,\halfrootthree)*={} **@{-},
"0" ; (-0.5,-\halfrootthree)*={} **@{-},
"0" ; (1,0)*={}="1" **@{-},
"1" ; (1.5,\halfrootthree)*={} **@{-},
"1" ; (1.5,-\halfrootthree)*={} **@{-},
\end{xy}
$$

Let us take two lattices $M_1$ and $M_2$ with distinguished bases
$\{e^i_0,\dots,e^i_{m_1}\}$, for $i=1,\ 2$ and respective subcube
polytopes $\Delta_i$, each of them having the set of vertices $\cA_i$.
Let $\ell_i=(e^i_0)^*: M_i\ra \Z$ be the projection to the zeroth
coordinate and, by abuse of notation, by the same letter we will
denote its composition with the projection $M_1\times M_2\ra
M_i\ra\Z$. Now we can take fiber product of each of these objects,
relative over the projection,
e.g. $M_1\phantom{*}_{\ell_1}\!\!\!\times_{\ell_2} M_2 \subset
M_1\times M_2$ consists of pairs $(u_1,u_2)$ such that
$\ell_1(u_1)=\ell_2(u_2)$. In other words
$M_1\phantom{*}_{\ell_1}\!\!\!\times_{\ell_2} M_2 =
\ker(\ell_1-\ell_2)$ and
$\Delta_1\phantom{*}_{\ell_1}\!\!\!\times_{\ell_2}\Delta_2=
\left(\Delta_1\times\Delta_2\right)\cap\ker(\ell_1-\ell_2)$

\end{const}
\begin{lemma}\label{product-subcube-is-subcube}
In the above situation
$\Delta=\Delta_1\phantom{*}_{\ell_1}\!\!\!\times_{\ell_2}\Delta_2$ is
a subcube polytope in $M=M_1\phantom{*}_{\ell_1}\!\!\!\times_{\ell_2}
M_2$ with the set of vertices
$\cA=\cA_1\phantom{*}_{\ell_1}\!\!\times_{\ell_2}\cA_2$. In general,
if $\Delta_i\subset(M_i)_\R$ and $\ell_i: M_i\ra \Z$ are lattice
homomorphisms such that $\ell_i(\Delta_i)\subset [0,1]$ then the set
of vertices of
$\Delta=\Delta_1\phantom{*}_{\ell_1}\!\!\!\times_{\ell_2}\Delta_2$ is
the fiber product of the vertices of $\Delta_i$'s.
\end{lemma}

\begin{pf}
The only non-trivial thing is to show that all vertices of $\Delta$
are in the fiber product of vertices of $\Delta_1$ and $\Delta_2$.
Since
$\Delta=\left(\Delta_1\times\Delta_2\right)\cap\ker(\ell_1-\ell_2)$ is
a codimension 1 linear section of $\Delta_1\times\Delta_2$ its
vertices are either vertices of $\Delta_1\times\Delta_2$ (which is
what we want) or are obtained by intersecting the hyperplane
$\ker(\ell_1-\ell_2)$ with an edge of $\Delta_1\times\Delta_2$.  To
this end, let us take two pairs of vertices $(u_1^1,u_2^1),\
(u_1^2,u_2^2)$, where $u_j^i$ is a vertex of $\Delta_j$. Suppose that
for some $t\in (0,1)$ the point $u=t(u_1^1,u_2^1)+(1-t)(u_1^2,u_2^2)$
is in $\ker(\ell_1-\ell_2)$, that is, we have
$$
\ell_1(tu_1^1+(1-t)u_1^2)=\ell_2(tu_2^1+(1-t)u_2^2)
$$ and moreover $\ell_1(u_1^1)\ne\ell_2(u_2^1)$. Thus, we may assume
that $\ell_1(u_1^1)=0$ and $\ell_2(u_2^1)=1$.  Hence, because of the
above equality and since $t\ne 0, 1$, we get $\ell_1(u^2_1)=1$,
$\ell_2(u^2_2)=0$ and $t=1/2$. So
$u=t(u_1^1,u_2^1)+(1-t)(u_1^2,u_2^2)=
t(u_1^1,u_2^2)+(1-t)(u_1^2,u_2^1)$ and $(u_1^1,u_2^2),\ (u_1^2,u_2^1)$
are vertices of $\Delta$ so $u$ lies in the interior of an edge of
$\Delta$.
\end{pf}

\begin{ex}\label{ex-product-of-tetrahedra}
Let us consider two copies of a tetrahedron, as in example
\ref{ex-tetrahedron}. That is, for $i=1, 2$, in a lattice $M^i=\Z
e^i_0\oplus\Z e^i_1 \oplus \Z e^i_2$, we consider a tetrahedron
$\Delta^3_i$ spanned on vertices $0$, $e^i_0+e^i_1$, $e^i_1+e^i_2$ and
$e^i_2+e^i_0$. We take the projections $(e^i_0)^*: M^i\ra \Z$ and in
the fiber product
$$M=M_1\phantom{*}_{\ell_1}\!\!\!\times_{\ell_2}
M_2=\ker\left((e_0^2-e_0^1)^*\right)\subset M_1\times M_2$$ by $e_0$
we denote the element $e^1_0+e^2_0$. The resulting fiber product of
tetrahedra $$\Delta=(\Delta^3_1\times\Delta^3_2)\cap M_\R$$ has the
following vertices: $0, e^1_1+e^1_2, e^2_1+e^2_2, e^1_1+e^1_2 +
e^2_1+e^2_2, e_0+e^1_1+e^2_1, e_0+e^1_1+e^2_2, e_0+e^1_2+e^2_1,
e_0+e^1_2+e^2_2$. 
\end{ex}

\begin{prop}\label{polytope-of-graft-is-fiberproduct}
Let $(\T_1,\ell_1)$ and $(\T_2,\ell_2)$ be two pointed trees. Then
$$
\hM(\T_1\phantom{*}_{\ell_1}\!\!\!\vee_{\ell_2}\T_2)=
\hM(\T_1)\phantom{*}_{\ell_1}\!\!\!\times_{\ell_2}\hM(\T_2)
$$
$$
\Delta(\T_1\phantom{*}_{\ell_1}\!\!\!\vee_{\ell_2}\T_2)=
\Delta(\T_1)\phantom{*}_{\ell_1}\!\!\!\times_{\ell_2}\Delta(\T_2)
$$

\end{prop}

\begin{pf}
Let $M_1=M(\T_1)$, $M_2=M(\T_2)$, and similarly for $N$'s, $\hM$' and
$\hN$'s. We set
$\T=\T_1\phantom{*}_{\ell_1}\!\!\!\vee_{\ell_2}\T_2$. Then, by
construction \ref{const-fiber-product},
$M=M(\T)=M_1\phantom{*}_{\ell_1}\!\!\!\times_{\ell_2} M_2$ and
$\boxdot_M=\boxdot_{M_1}
\phantom{.}_{\ell_1}\!\!\!\times_{\ell_2}\boxdot_{M_2}$.  The two
projections $p_i: M\ra M_i$ yield respective injections of $\Hom(\
\cdot\ ,\Z)$-spaces: $\iota_i: N_i\hookrightarrow N$ (in fact
$N=(N_1\times N_2)/\Z(\ell_1-\ell_2)$). If $\N_i$ and $\N$ denote,
respectively, inner nodes of $\T_i$ and $\T$ then
$\N=\iota_1(\N_1)\cup \iota_2(\N_2)$.  Since $\hN$, $\hN_1$, $\hN_2$
are defined by extending $N$, $N_1$, $N_2$ by $\N/2$, $\N_1/2$ and
$\N_2/2$, respectively, it follows that $\hN=\hN_1+\hN_2$ in
$N_\R$. This implies the first equality of the lemma. Similarly, since
the set $\N$ determines vertices of $\boxdot_M$ which span
$\Delta(\T)$, see \ref{polytope-model}, we get the second equality.
\end{pf}

The above result can be expressed as follows: the polygon of a tree
$\T$ is a fiber product of polygons of star trees associated to inner
nodes of the tree, fibered over the relations encoded in the inner
branches of the tree. Since $\Delta(\T^{\!\!*3})$ is a 3-dimensional
tetrahedron this is especially straightforward in case of 3-valent
trees.

For any inner node $v\in\N$ of a 3-valent tree $\T$ we consider the
lattice $M_v=\Z e^v_1\oplus\Z e^v_2\oplus\Z e^v_3$, where $e_1$, $e_2$
and $e_3$ are the three edges stemming from $v$. Inside $M_v$ we have
the tetrahedron $\Delta_v$ with vertices $0$, $e^v_1+e^v_2$,
$e_2^v+e_3^v$ and $e_3^v+e_1^v$. We consider the big lattice
$\widetilde{M}=\bigoplus_{v\in\N}M_v$ and $\widetilde{M}_\R$ contains
the product $\prod_{v\in\N}\Delta_v$. Now for each inner edge $e$ we
have a form $\tilde{e}^*: \widetilde M\ra \Z$ such that
$\tilde{e}^*(e^v_i)=(-1)^\epsilon$ if $\partial_\epsilon(e)=v$ and
$e^v_i=e$, and $\tilde{e}^*(e^v_i)=0$ otherwise. Then the intersection
$\bigcap_{e\in\E^o}\ker(\tilde{e}^*)$ can be identified with the
lattice $M$, that is, we map $e\in\E^o$ to $e^{v}_i+e^{v'}_{i'}$ where
$\partial(e)=(v,v')$, and $e^{v}_i=e^{v'}_{i'}=e$, while a petiole $e$
is just mapped to its unique representation in $\widetilde M$. Then,
by \ref{polytope-of-graft-is-fiberproduct} we get
$$\Delta(\T)=\left(\prod_{v\in\N}\Delta_v\right)\cap
\left(\bigcap_{e\in\E^o}\ker(\tilde{e}^*)_\R\right)$$


\subsection{Geometric model of a tree}\label{const-toric-var}

First, let us recall the construction of a projective tori variety
from a lattice polytope of characters. Let $\widehat{M}$ and $\hN$ be
dual lattices of characters and 1 parameter subgroups for a torus
$T_\hN=\hN\otimes_\Z\C^*$.

\begin{defin}\label{normal-polytope}
A lattice polytope $\Delta\subset \widehat{M}_\R$ is called normal if
\begin{itemize}
\item the sublatice of $\widehat{M}$ spanned by the differences of
points in $\Delta\cap\widehat{M}$ is equal to $\widehat{M}$
\item for every integer $d\geq 0$ any point in $d\Delta\cap\widehat{M}$ is
equal to a sum of $d$ points in $\Delta\cap\widehat{M}$.
\end{itemize}

\end{defin}

Equivalently, the second condition in the above definition can be
restated as follows. Let $\widehat{M}'=\widehat{M}\oplus \Z$ and take
an affine map $i_1: \widehat{M}\ra\widehat{M}'$ such that
$i_1(u)=(u,1)$. Then $\Delta$ is normal in $\hM$ if and only if the
semigroup spanned in $\widehat{M}'$ by $i_1(\Delta\cap\widehat{M})$ is
equal to the semigroup of lattice points in cone spanned in $\hM_\R'$
by $i_1(\Delta)$, that is the semigroup $\R_{\geq
0}(i_1(\Delta))\cap\widehat{M}'$.

\begin{defin}\label{def-toric-via-polytope}
Suppose that $\Delta$ is a normal polytope in $\hM$. Let $A^d_\Delta$
be a $\C$-linear space with the basis $\{\chi^u: u\in
d\Delta\cap\widehat{M}\}$. We consider a graded $\C$-algebra
$A(\Delta)=\bigoplus_{d\geq 0}A^d_\Delta$, with multiplication
$\chi^{u_1}\chi^{u_2}=\chi^{u_1+u_2}$. Then $X(\Delta)=\Proj
A(\Delta)$ is called the projective model of $\Delta$.
\end{defin}

We note that in the above situation $A(\Delta)$ is a normal ring, that
is, it integrally closed in its field of fractions.  This, by
definition, is equivalent to saying that affine spectrum
$\Spec(A(\Delta))$ is a normal affine variety. In fact, in such a case
$A(\Delta)$ is the semigroup algebra of $\R_{\geq
0}(i_1(\Delta))\cap\widehat{M}'$ so $\Spec(A(\Delta))$ is an affine
toric variety with the big torus $T_{\hN\oplus\Z}$. In the projective
case we have the following general result which summarizes properties
of the projective model of a normal polytope, see
\cite[Sect.~2.1--2.4]{oda}, \cite{sturmfels} or \cite{fulton}.

\begin{prop}\label{toric-via-polytope}
Suppose that $\Delta$ is a normal polytope in the lattice
$\widehat{M}$ of characters of a torus $T_{\widehat{N}}$.  
Then the following holds:
\begin{enumerate}
\item $X(\Delta)$ is a toric variety on which $T_{\hN}$ acts effectively,
\item $X(\Delta)$ is embedded in $\P^{|\Delta\cap\widehat{M}|-1}$ as a
projectively normal variety such that $H^0(X(\Delta),\cO_X(d))=A^d_\Delta$,
\item Characters from $\Delta\cap\widehat{M}$ define a diagonal action
of $T_{\widehat{N}}$ on $\P^{|\Delta\cap\widehat{M}|-1}$ which
restricts to the torus action on $X(\Delta)$,
\item The induced action of $T_{\widehat{M}}$ on
$H^0(X(\Delta),\cO_X(d))$ is linearizable with weights in
$d\Delta\cap\widehat{M}$.
\item $X(\Delta)\subset\P^{|\Delta\cap\widehat{M}|-1}$ is the closure of the
image of the map $T_{\widehat{M}}\ra\P^{|\Delta\cap\widehat{M}|-1}$
defined by the characters from $\Delta\cap\widehat{M}$
\end{enumerate}
\end{prop}

Because of \ref{tree-polytope-normal} the polytope model $\Delta(\T)$
of a {\bf 3-valent} tree $\T$ is normal so we can consider its
projective model. The following is the key result of the paper which
allows us to study projective models of binary symmetric trees in
purely toric way.

\begin{thm}\label{equivalence-of-models}
Let $(\T,W,\hW)$ be a binary symmetric {\bf 3-valent} tree.
Then the varieties $X(\T)$ and $X(\Delta(\T))$ are projectively equivalent
in $\P(W_\L^\rho)=\P^{2^{|\L|-1}-1}$.
\end{thm}

\begin{pf}
By \ref{polytope=>parametrization} the parametrization of $X(\T)$ is
defined as a rational map from $\product_{e\in\E}\P(\hW^e)$ defined by
characters of torus $T_N$ which are vertices of $\Delta(\T)$. Thus,
$X(\T)$ is the closure of the respective map $T_N\ra
\P^{2^{|\L|-1}-1}$. Since $\hM\subset M$ is the sublattice spanned by
vertices of $\Delta(\T)$ this factors to the map $T_\hN\ra
\P^{2^{|\L|-1}-1}$ the image of which defines $X(\Delta)$,
\ref{toric-via-polytope}.
\end{pf}


\subsection{1-parameter group action, quotients.}\label{sect-quotients}


In this section we consider quotients of projective varieties as in
Mumford's GIT \cite{git}. For a comprehensive exposition of the
theory, including a relevant definition of good quotient we refer to
\cite{abb}. In the present section as well as in section
\ref{sect-mutation-deformation} we consider an algebraic action of a
torus $T$ on a projective variety $X\hookrightarrow \P^m$ which is
given by a choice of weights hence it extends to the affine cone over
$X$ and thus it determines its linearization, its set of semi-stable
points $X^{ss}$ and its good quotient $X^{ss}\ra X^{ss}// T$, see
\cite[Ch.6]{abb}.

\begin{const}\label{construction-graft-as-quotient} 

Let $\Delta_i\subset (\hM_i)_\R$, for $i=1, 2$ be two lattice
polytopes admitting unimodular covers hence normal, see
\ref{sect-unimodular}, and $X(\Delta_i)\subset \P^{n_i-1}$, where
$n_i=|\hM_i\cap\Delta_i|$, their associated toric varieties.  In
$M^\times=\hM_1\times\hM_2$ we take the product polytope
$\Delta^\times=\Delta_1\times\Delta_2$ which is also normal,
ref{product-unimodular}. Then the associated toric variety
$X^\times=X(\Delta^\times)\subset \P^{n_1n_2-1}$ is the Segre image of
$X(\Delta_1)\times X(\Delta_2)$.

Suppose that $\ell_i: \hM_i \ra \Z$ are lattice homomorphisms such
that $(\ell_i)_\R(\Delta_i)\subset [0,1]$.  We pull $\ell_i$ to the
product of lattices and on $\hM_1\times\hM_2$ we define the form
$(\ell_1-\ell_2)$. The form defines a diagonal action
$\lambda^{\ell_1-\ell_2}$ of $\C^*$ on $X^\times\subset\P^{n_1n_2-1}$
which on the coordinate associated to $\chi^{(u_1,u_2)}$, where
$u_i\in\Delta_i\cap\hM_i$, has the weight
$\ell_1(u_1)-\ell_2(u_2)\in\{-1,0,1\}$. Accordingly, we regroup the
coordinates of $\P^{n_1n_2-1}$ and write them as $[z^-_i,z^0_j,z^+_k]$
depending on whether they are of weight $-1$, $0$ and $1$,
respectively. That is
$$\lambda_{\ell_1-\ell_2}(t)[z^-_i,z^0_j,z^+_k]= [t^{-1}z^-_i,z^0_j,t
z^+_k]$$ The above formula defines the action of
$\lambda_{\ell_1-\ell_2}$ on the cone over $X^\times$ and thus a
$\C^*$-linearization of the bundle $\cO_{X^\times}(1)$ in the sense of
GIT.  By $X^0$ let us denote the intersection of $X^\times$ with the
complement of the space spanned on the eigenvectors of
$\lambda_{\ell_1-\ell_2}$ of weight $\ne 0$, that is
$X^0=X^\times\setminus\{[z_i^-,z^0_j,z_k^+]: \forall j\ z^0_j=0\}$

We set $\hM=\ker(\ell_1-\ell_2)$ and $\Delta=\Delta^\times\cap
\ker(\ell_1-\ell_2)= \Delta_1\phantom{*}_{\ell_1}\!\!\!\times_{\ell_2}
\Delta_2$. By \ref{fiber-product-unimodular-gen} $\Delta$ is a normal
polytope and by $X(\Delta)$ we denote its associated toric variety.
\end{const}

\begin{prop}\label{thm-graft-as-quotient}
In the above situation the set $X^0$ is equal to the set of the
semistable points of the action of $\lambda^{\ell_1-\ell_2}$.  The
projection to the weight 0 eigenspace $[z^0_i,z^0_j,z^+_k]\mapsto
[z^0_j]$ defines a regular map of $X^0$ to $X(\Delta)$ and $X(\Delta)$
is a good quotient for the action of $\lambda_{\ell_1-\ell_2}$.
\end{prop}

\begin{pf}
The sections of $\cO_{X^\times}(m)$ for $m>0$ make a vector space
spanned on $\chi^u$, where $u\in m\Delta^\times\cap M^\times$. Among
them, these which are invariant with respect to the action of
$\lambda^{\ell_1-\ell_2}$ are associated to $u$'s in the intersection
with $\ker{\ell_1-\ell_2}$ thus in $m\Delta\cap M$. By the normality
of $\Delta$, see \ref{fiber-product-unimodular-gen}, the algebra of
invariant sections is generated by these from
$\cO_{X^\times}(1)$. Thus the set of semistable points of the action
of $\lambda_{\ell_1-\ell_2}$ is where at least one of the coordinates
$z^0_j$ is non-zero and the quotient map is the projection to the
weight zero eigenspace.
\end{pf}

\begin{cor}
Let $(\T_1,\ell_1)$ and $(\T_2,\ell_2)$ be two pointed trees. Then
$X(\T_1\phantom{*}_{\ell_1}\!\!\!\vee_{\ell_2}\T_2)$ is a good
quotient of $X(\T_1)\times X(\T_2)$ with respect to an action of
$\lambda^{\ell_1-\ell_2}$.
\end{cor}

\begin{ex}\label{ex-quotient-of-two-P2}
Consider the $\C^*$ action on the product $\P^3_1\times\P^3_2$
given by the formula:
$$\lambda(t) ([z_0^1,z_1^1,z_2^1,z_3^1],[z_0^2,z_1^2,z_2^2,z_3^2]) =
([z_0^1,t z_1^1,t z_2^1,z_3^1],[z_0^2, t^{-1}z_1^2,
t^{-1}z_2^2, z_3^2])$$ where the superscripts of the coordinates
indicate the factor in the product $\P^3_1\times\P^3_2$. The following
rational map $\P^3_1\times\P^3_2 \ - \to \P^7$ is $\lambda$
equivariant and regular outside the set $\{
z_0^1=z_3^1=z_1^2=z_2^2=0\}\cup\{ z_1^1=z_1^2=z_0^2=z_3^2=0\}$, each
component of this set is a quadric $\P^1\times\P^1$:
$$([z_0^1,z_1^1,z_2^1,z_3^1],[z_0^2,z_1^2,z_2^2,z_3^2]) \mapsto
[z_0^1z_0^2,z_0^1z_3^2,z_1^1z_1^2,z_1^1z_2^2,z_2^1z_1^2,
z_2^1z_2^2,z_3^1z_0^2,z_3^1z_3^2]$$ If $[x_0,\dots x_7]$ are
coordinates in $\P^7$ then the image of this map is the intersection
of two quadrics $\{ x_0x_7=x_1x_6\}\cap\{ x_2x_5=x_3x_4\}$.

The above claim will be clear if we write functions $z_i^1z_j^2$ in
terms of characters of the respective torus, which we denote by
$e^1_i$ and $e^2_j$, respectively.  Namely, dividing the right hand
side of the above displayed formula by $z_0^1z_0^2$ we get the following
sequence of rational functions:
$$
[1, \chi^{e_1^2+e_2^2}, \chi^{e_0^1+e_1^1}\chi^{e_0^2+e_1^2}, 
\chi^{e_0^1+e_1^1}\chi^{e_0^2+e_2^2}, \chi^{e_0^1+e_2^1}\chi^{e_0^2+e_1^2}, 
\chi^{e_0^1+e_2^1}\chi^{e_0^2+e_2^2}, \chi^{e_1^1+e_2^1}, 
\chi^{e_1^1+e_2^1}\chi^{e_1^2+e_2^2}]
$$ If we write the sums of the exponents of the above rational
functions in $M_1\oplus M_2$ and call $e_0=e^1_0+e^2_0$ then we get
the vertices of $\Delta(\T^{\!\!*3}\vee\T^{\!\!*3})$ which we computed
in example \ref{ex-product-of-tetrahedra}. From the above formula we
can read the weights with which 1-parameter groups
$\lambda_{(e_j^i)^*}$, for $i,\ j=1,\ 2$, associated to leaves, act on
the quotient variety in $\P^7$.

\end{ex}


\section{3-valent binary trees.}


{\bf Summary:} From this point on we concentrate on understanding
varieties associated to 3-valent binary trees and we prove main
results of the present note which are as follows: (1) such varieties
have only Gorenstein terminal singularities and are Fano of index 4,
(2) any two such varieties associated to trees with the same number of
leaves are in the same connected component of the Hilbert scheme of
the projective space, (3) their Hilbert-Ehrhard polynomial can be
computed effectively..

\subsection{Paths, networks and sockets.}\label{sect-toric-ideal}

Let $\T$ be a 3-valent binary symmetric tree.  In section
\ref{sect-polytope-model} we identified the variety $X(\T)$ in
$\P^{2^{|\L|-1}-1}$ with the closure of the image of a torus map
defined by a polytope $\Delta(\T)$. We recall that the linear
coordinates on the ambient projective space can be identified with the
vertices of $\Delta(\T)$ which are among the vertices of the cube
$\boxdot_M$ satisfying parity relation with respect to the forms
$v\in\N\subset N$, \ref{polytope-model}. For 3-valent trees we have a
convenient interpretation of these points.

\begin{defin}\label{def-path-netowork-socket}
Let $\T$ be a 3-valent tree. A {\em path} $\gamma$ on $\T$ of length $m\geq
1$ is a choice of $m+1$ distinct vertices $v_0,\dots,v_m$ such that
$v_0$ and $v_m$ are leaves (called the ending points of $\gamma$) and
there exists $m$ edges, $e_1,\dots,e_m$ such that for $i=1,\dots,m$ it
holds $\partial(e_i)=\{v_{i-1},v_i\}$. 

A {\em network} of paths (or just a network) $\Gamma$ on $\T$ is a set
of paths (possibly an empty set), each two of them have no common
vertex (neither edge). For any network of paths $\Gamma$ on $\T$ we
define the {\em socket} $\mu(\Gamma)\subset\cL$ to be the set of
leaves which are ending points of paths in $\Gamma$.

A tree $\T$ is {\em labeled}  if its leaves are
numbered by $1,\dots, |\cL|$.  A subset $\mu\subset\cL$ is represented
by a characteristic sequence $\kappa(1),\dots,\kappa(|\cL|)$ in which
$\kappa(i)=1$ or $0$, depending on whether the leaf numbered by $i$ is
in $\mu$ or not.
\end{defin}

Sockets of networks will identified by their characteristic binary
sequences. We note that, clearly, every socket consists of even number
of elements in $\cL$.

\begin{ex}\label{4-leaf-tree-networks}
Let us consider a labeled 3-valent tree with four leaves. In the
following diagram, in the upper row we draw all possible networks on
this tree, where paths are denoted by solid line segments. In the
lower row we write down the respective sockets in terms of
characteristic sequences of length four

$$
\begin{array}{cccccccc}
\begin{xy}<7pt,0pt>:
(0,0)*={}="0" ; (-1,\rootthree)*={1} **@{.},
"0" ; (-1,-\rootthree)*={2} **@{.},
"0" ; (2,0)*={}="1" **@{.},
"1" ; (3,\rootthree)*={3} **@{.},
"1" ; (3,-\rootthree)*={4} **@{.},
\end{xy}
&
\begin{xy}<7pt,0pt>:
(0,0)*={}="0" ; (-1,\rootthree)*={1} **@{-},
"0" ; (-1,-\rootthree)*={2} **@{-},
"0" ; (2,0)*={}="1" **@{.},
"1" ; (3,\rootthree)*={3} **@{.},
"1" ; (3,-\rootthree)*={4} **@{.},
\end{xy}
&
\begin{xy}<7pt,0pt>:
(0,0)*={}="0" ; (-1,\rootthree)*={1} **@{.},
"0" ; (-1,-\rootthree)*={2} **@{.},
"0" ; (2,0)*={}="1" **@{.},
"1" ; (3,\rootthree)*={3} **@{-},
"1" ; (3,-\rootthree)*={4} **@{-},
\end{xy}
&
\begin{xy}<7pt,0pt>:
(0,0)*={}="0" ; (-1,\rootthree)*={1} **@{-},
"0" ; (-1,-\rootthree)*={2} **@{-},
"0" ; (2,0)*={}="1" **@{.},
"1" ; (3,\rootthree)*={3} **@{-},
"1" ; (3,-\rootthree)*={4} **@{-},
\end{xy}
&
\begin{xy}<7pt,0pt>:
(0,0)*={}="0" ; (-1,\rootthree)*={1} **@{-},
"0" ; (-1,-\rootthree)*={2} **@{.},
"0" ; (2,0)*={}="1" **@{-},
"1" ; (3,\rootthree)*={3} **@{-},
"1" ; (3,-\rootthree)*={4} **@{.},
\end{xy}
&
\begin{xy}<7pt,0pt>:
(0,0)*={}="0" ; (-1,\rootthree)*={1} **@{-},
"0" ; (-1,-\rootthree)*={2} **@{.},
"0" ; (2,0)*={}="1" **@{-},
"1" ; (3,\rootthree)*={3} **@{.},
"1" ; (3,-\rootthree)*={4} **@{-},
\end{xy}
&
\begin{xy}<7pt,0pt>:
(0,0)*={}="0" ; (-1,\rootthree)*={1} **@{.},
"0" ; (-1,-\rootthree)*={2} **@{-},
"0" ; (2,0)*={}="1" **@{-},
"1" ; (3,\rootthree)*={3} **@{-},
"1" ; (3,-\rootthree)*={4} **@{.},
\end{xy}
&
\begin{xy}<7pt,0pt>:
(0,0)*={}="0" ; (-1,\rootthree)*={1} **@{.},
"0" ; (-1,-\rootthree)*={2} **@{-},
"0" ; (2,0)*={}="1" **@{-},
"1" ; (3,\rootthree)*={3} **@{.},
"1" ; (3,-\rootthree)*={4} **@{-},
\end{xy}
\\
\\
0,0,0,0&1,1,0,0&0,0,1,1&1,1,1,1&1,0,1,0
&1,0,0,1&0,1,1,0&0,1,0,1
\end{array}
$$
\end{ex}

\begin{lemma}\label{networks=vertices}
Let $\T$ be a 3-valent tree. Associating to a network $\Gamma$ a point
$u(\Gamma)=\sum_e \Gamma(e)\cdot e\in M(\T)$, where $\Gamma(e)=1, 0$
depending on whether $e$ is on $\Gamma$ or not, defines a bijection
between networks and vertices of $\Delta(\T)$
\end{lemma}

\begin{pf}

First note that $u(\Gamma)\in\Delta(\T)$.  To define the inverse of
$\Gamma\ra u(\Gamma)$, for any vertex $u=\sum_{e\in \E}
\epsilon_e\cdot e\in\Delta(\T)$ we define the support of $u$
consisting of edges of $\T$ whose contribution to $u$ is nonzero,
i.e. $\{e\in\E : e^*(u)=1\}$.  The parity condition $\forall_{v\in\N}$
either $v(u)=0$ or $v(u)=2$ yields that these edges define a network
on $\T$.
\end{pf}

We note that, because of \ref{number-of-vertices}, there are
$2^{|\L|-1}$ networks. On the other hand, the association of the
socket to a network gives a map from the set of networks to the
subsets of leaves.  This map is surjective, that is, every subset
$\mu$ of $\L$ with even umber of elements is a socket of a
network. Indeed, this follows by a straightforward induction with
respect to the number of leaves of the tree: in the induction step we
write a tree $\T_{n+1}$ with $n+1$ leaves as a graft of a tree $\T_n$
with $n$ leaves and a star tree with 3 leaves and consider three cases
depending on how many of the two new leaves replacing one old are in
the set $\mu\subset\L$.

Finally, because the number of all subsets of $\L$ with even number of
elements equals to $2^{|\L|-1}$ we get the following.

\begin{lemma}\label{vertices=sockets}
Let $\T$ be a 3-valent tree. Then associating to a network its
socket defines a bijection between the set of networks of paths on
$\T$ and the set of subsets of $\L$ which have even number of
elements.
\end{lemma}

We note that the sockets of a tree $\T$ form a convenient basis in the
space $W_\cL^\rho$, which was introduced in section
\ref{bin-symm-trees}.  Indeed, in order to use toric arguments we have
diagonalized the action of the involution $\rho$ on $W$ with a basis
$\nu_0$, $\nu_1$ such that $\rho(\nu_i)=(-1)^i\nu_i$. Now any socket
(or, equivalently, a subset of $\L$ with even number of elements)
whose characteristic binary function is $\kappa: \cL\ra\{0,\ 1\}$,
defines an element $\otimes_{v\in\cL} \nu_{\kappa(v)}^*$ in
$(W_\cL^\rho)^*$. Similarly, to any network$\Gamma$ on $\T$ we
associate a vector $\otimes_{e\in\E} \omega_{\Gamma(e)}^*$ in
$\widehat{W}^\E$, where $\omega_i$ is such
$\widehat\rho(\omega_i)=(-1)^i\omega_i$ and $\Gamma(e)=1, 0$ depending
on whether $e$ is in $\Gamma$ or not. Now associating to a network its
sockets defines an isomorphism $(\widehat{W}^\E)^{G_\N}\ra
(W_\cL^\rho)$ which one can compare to what we discuss in \ref{map-injective}.

We have a convenient description of the action of one-parameter groups
associated to leaves of $\T$ in terms of socket coordinates of
$\P(W_\cL^\rho)$ . Namely, given a leaf $\ell$ the 1-parameter group
$\lambda^\ell$ acts on the coordinate $\chi^\kappa$ with the weight
$\kappa(\ell)$.

\begin{const}\label{toric-intersection}
Using networks and sockets, and the toric formalism, one can explain
the inclusion $X(\T)\subset \P(W_\cL^\rho)$ as follows. Let
$\widetilde{M}=\bigoplus_{\kappa\ne 0}\Z \cdot\kappa$ be a lattice, a
free abelian group generated by non-empty sockets of a tree $\T$.  The
empty socket $\kappa=0$ we interpret as the zero of the lattice.  Then
$\P((W_\cL^\rho)$ is a toric variety $X(\widetilde\Delta^0)$
associated to a unit simplex $\widetilde\Delta^0$ in $\widetilde M$
spanned on the vectors of the distinguished basis. 

Now the bijective map {\em sockets} $\leftrightarrow$ {\em networks}
gives rise to a homomorphism of lattices
$\widetilde{M}\ra\widehat{M}$, where, recall, the latter lattice is
spanned in $M$ by the points associated to networks. This gives a
surjective map from the symmetric graded algebra spanned by all the
sockets, which is just algebra of polynomials $\C[\chi^\kappa]$, to
the algebra $A(\Delta)$, hence we get the inclusion $X(\T)\subset
\P(W_\cL^\rho)$, c.f.~\ref{const-toric-var} and
\ref{toric-via-polytope}.
\end{const}

\begin{defin}\label{def-primitive-rel}
Let $\Delta$ be a normal lattice polytope in a lattice $M$. Let us
choose two collection of points $u_1,\dots,u_r$ and $w_1,\dots,w_s$ in
$\Delta\cap M$ and positive integers $a_1,\dots,a_r$ and
$b_1,\dots,b_s$. This data defines a {\em relation} of degree $d$ for
$\Delta$ if $a_1+\cdots+a_r=b_1+\cdots+b_s=d$ and
$$a_1u_1+\cdots+a_ru_r=b_1w_1+\cdots+b_sw_s$$ 
The relation is called
primitive if $\{u_1,\dots,u_r\}\cap\{w_1,\dots,w_s\}=\emptyset$.
\end{defin}

Let us recall that given the projective variety $X\subset \P^r$ with
graded coordinate ring $S(X)=\bigoplus_{m\geq 0} S^m(X)$ its ideal
$\cI(X)$ is the kernel of evaluation map $Symm(S^1(X))\ra S(X)$.  The
following result is known as binomial generation of a toric ideal, see
\cite{es}, \cite{sturmfels}.

\begin{lemma}\label{toric-ideal-via-relations}
Suppose that we are in the situation of \ref{toric-via-polytope}.
Then the ideal $\cI(X(\Delta))$ is generated by polynomials
$$
\left(\chi^{u_1}\right)^{a_1}\cdots\left(\chi^{u_r}\right)^{a_r}-
\left(\chi^{w_1}\right)^{b_1}\cdots\left(\chi^{w_s}\right)^{b_s}
$$ 
where $u_1,\dots,u_r$ and $w_1,\dots,w_s$, together with
$a_1,\dots,a_r$ and $b_1,\dots,b_s$ define a primitive relation
for $\Delta$.
\end{lemma}

\begin{ex}\label{4-leaf-tree-relations}
 The following are primitive relations and respective equations
for the polytope coming from a 3-valent tree with four leaves,
c.f.~example \ref{4-leaf-tree-networks}.  First, we describe them in
terms of networks; they are as follows:
$$
\begin{xy}<10pt,0pt>:
(0,0)*={}="0" ; (-0.5,\halfrootthree)*={1} **@{.},
"0" ; (-0.5,-\halfrootthree)*={2} **@{.},
"0" ; (1,0)*={}="1" **@{.},
"1" ; (1.5,\halfrootthree)*={3} **@{.},
"1" ; (1.5,-\halfrootthree)*={4} **@{.},
\end{xy}
\ \ + \ \ 
\begin{xy}<10pt,0pt>:
(0,0)*={}="0" ; (-0.5,\halfrootthree)*={1} **@{-},
"0" ; (-0.5,-\halfrootthree)*={2} **@{-},
"0" ; (1,0)*={}="1" **@{.},
"1" ; (1.5,\halfrootthree)*={3} **@{-},
"1" ; (1.5,-\halfrootthree)*={4} **@{-},
\end{xy}
\ \ = \ \ 
\begin{xy}<10pt,0pt>:
(0,0)*={}="0" ; (-0.5,\halfrootthree)*={1} **@{-},
"0" ; (-0.5,-\halfrootthree)*={2} **@{-},
"0" ; (1,0)*={}="1" **@{.},
"1" ; (1.5,\halfrootthree)*={3} **@{.},
"1" ; (1.5,-\halfrootthree)*={4} **@{.},
\end{xy}
\ \ + \ \ 
\begin{xy}<10pt,0pt>:
(0,0)*={}="0" ; (-0.5,\halfrootthree)*={1} **@{.},
"0" ; (-0.5,-\halfrootthree)*={2} **@{.},
"0" ; (1,0)*={}="1" **@{.},
"1" ; (1.5,\halfrootthree)*={3} **@{-},
"1" ; (1.5,-\halfrootthree)*={4} **@{-},
\end{xy}
$$
and
$$
\begin{xy}<10pt,0pt>:
(0,0)*={}="0" ; (-0.5,\halfrootthree)*={1} **@{-},
"0" ; (-0.5,-\halfrootthree)*={2} **@{.},
"0" ; (1,0)*={}="1" **@{-},
"1" ; (1.5,\halfrootthree)*={3} **@{.},
"1" ; (1.5,-\halfrootthree)*={4} **@{-},
\end{xy}
\ \ + \ \
\begin{xy}<10pt,0pt>:
(0,0)*={}="0" ; (-0.5,\halfrootthree)*={1} **@{.},
"0" ; (-0.5,-\halfrootthree)*={2} **@{-},
"0" ; (1,0)*={}="1" **@{-},
"1" ; (1.5,\halfrootthree)*={3} **@{-},
"1" ; (1.5,-\halfrootthree)*={4} **@{.},
\end{xy}
\ \ = \ \ 
\begin{xy}<10pt,0pt>:
(0,0)*={}="0" ; (-0.5,\halfrootthree)*={1} **@{-},
"0" ; (-0.5,-\halfrootthree)*={2} **@{.},
"0" ; (1,0)*={}="1" **@{-},
"1" ; (1.5,\halfrootthree)*={3} **@{-},
"1" ; (1.5,-\halfrootthree)*={4} **@{.},
\end{xy}
\ \ + \ \ 
\begin{xy}<10pt,0pt>:
(0,0)*={}="0" ; (-0.5,\halfrootthree)*={1} **@{.},
"0" ; (-0.5,-\halfrootthree)*={2} **@{-},
"0" ; (1,0)*={}="1" **@{-},
"1" ; (1.5,\halfrootthree)*={3} **@{.},
"1" ; (1.5,-\halfrootthree)*={4} **@{-},
\end{xy}
$$
On $\P(W_\cL^\rho)$ we introduce coordinates
$x_{\kappa(1)\cdots \kappa(4)}$ indexed by characteristic sequences for the
sockets in $\cL$. Then the respective equations defining $X(\T)$ 
are as follows:
$$
\begin{array}{ccc}
x_{0000}\cdot x_{1111}&=&x_{1100}\cdot x_{0011}\\
x_{1001}\cdot x_{0110}&=&x_{1010}\cdot x_{0101}
\end{array}
$$ 
Finally, let us note that renumbering the leaves or, equivalently,
changing the shape of a 3-valent tree connecting the four numbered
leaf vertices, produces the following respective equations
$$
\begin{array}{lccr}
\begin{xy}<10pt,0pt>:
(0,0)*={}="0" ; (-0.5,\halfrootthree)*={1} **@{-},
"0" ; (-0.5,-\halfrootthree)*={3} **@{-},
"0" ; (1,0)*={}="1" **@{-},
"1" ; (1.5,\halfrootthree)*={2} **@{-},
"1" ; (1.5,-\halfrootthree)*={4} **@{-},
\end{xy}
& & &
\begin{array}{ccc}
x_{0000}\cdot x_{1111}&=&x_{1010}\cdot x_{0101}\\
x_{1001}\cdot x_{0110}&=&x_{1100}\cdot x_{0011}
\end{array}
\\ \\
\begin{xy}<10pt,0pt>:
(0,0)*={}="0" ; (-0.5,\halfrootthree)*={1} **@{-},
"0" ; (-0.5,-\halfrootthree)*={4} **@{-},
"0" ; (1,0)*={}="1" **@{-},
"1" ; (1.5,\halfrootthree)*={2} **@{-},
"1" ; (1.5,-\halfrootthree)*={3} **@{-},
\end{xy}
& & &
\begin{array}{ccc}
x_{0000}\cdot x_{1111}&=&x_{1001}\cdot x_{0110}\\
x_{1001}\cdot x_{0110}&=&x_{1010}\cdot x_{0101}
\end{array}
\end{array}
$$ We note that all the above equations involve only four quadratic
monomials: $x_{0000}x_{1111}$, $x_{1100}x_{0011}$, $x_{1010}x_{0101}$,
$x_{1001}x_{0110}$. Moreover, given any leaf $\ell$, the 1-parameter
group $\lambda^\ell$ acts with weight 1 on each of them.

\end{ex}


\subsection{Dual polytopes, fans, resolutions and Fano varieties.}


In the situation of \ref{def-toric-via-polytope} the description of
the fan of the variety $X(\Delta)$ in $\hN$ is given in terms of its
support functions \cite[Thm.~2.22]{oda} or dual polytopes
\cite{fulton}.

\begin{ex}\label{ex-tetrahedron1}
By looking at the example \ref{ex-tetrahedron} and the inequalities
which appear there we see that the fan of $X(\T\!\!^{*3})$ in $\widehat
N\supset N$ has rays generated by the following elements:
$-v/2=-(e^*_0+e^*_1+e_2^*)/2$, $v/2-e_0^*=(e_1^*+e_2^*-e_0^*)/2$,
$v/2-e_1^*=(e_0^*+e_2^*-e_1^*)/2$, $v/2-e_2^*=(e_0^*+e_1^*-e_2^*)/2$.
\end{ex}

The formula from \ref{polytope-of-graft-is-fiberproduct} can be used
to get the description of the polytope dual to $\Delta(\T)$, hence to
describing the fan of $X(\T)$ for 3-valent trees.

\begin{lemma}\label{defining-equations}
Let $\T$ be a 3-valent binary symmetric tree with $n$ inner nodes.
Then the polytope $\Delta(\T)$ is defined in $M_\R$ by $4n$
inequalities, which are as follows: for any inner node $v\in\N$, such
that $v=e_{0.1}^*+e_{v.1}^*+e_{v.2}^*$ we take
$$ (-v/2)(\ \cdot \ )\geq -1, \ \ (v/2 - e_{v.0}^*)(\ \cdot \ )\geq 0,
\ \ (v/2 - e_{v.1}^*)(\ \cdot \ )\geq 0, \ \ (v/2 - e_{v.2}^*)(\ \cdot
\ )\geq 0.
$$ 
\end{lemma}

\begin{pf}
Let $(\T_1,\ell_1)$ and $(\T_2,\ell_2)$ be pointed trees.  If
$\Delta_i=\Delta(\T_i)\subset (M_i)_\R$ is defined by inequalities
with respect to some forms $w^i_j$ in $(N_i)_\R$ then
$\Delta_1\times\Delta_2$ is defined by forms $(w^1_j,0)$ and
$(0,w_j^2)$ in $(N_1)_\R\times(N_2)_\R$.  Then the classes of these
forms in $N=(N_1\times N_2)_\R/\R(\ell_1-\ell_2)$ define
the fiber product of $\Delta_i$'s.
\end{pf}

\begin{defin}\label{def-dual-polytope}
For a binary symmetric 3-valent tree $\T$ we define a polytope
$\Delta^\vee(\T)$ in $N_\R$ which is the convex hull of
$-v/2=-(e^*_{v.0}+e^*_{v.1}+e_{v.2}^*)/2$,
$v/2-e_{v.0}^*=(e_{v.1}^*+e_{v.2}^*-e_{v.0}^*)/2$,
$v/2-e_{v.1}^*=(e_{v.0}^*+e_{v.2}^*-e_{v.1}^*)/2$,
$v/2-e_{v.2}^*=(e_{v.0}^*+e_{v.1}^*-e_{v.2}^*)/2$, for $v\in\N$ and
$e_{v.0}$, $e_{v.1}$, $e_{v.2}$ edges containing $v$.
\end{defin}

Let us note that the listed above points are in fact vertices of
$\Delta^\vee(\T)$. Indeed, take $v\in\N$ and $e_{v.0}$, $e_{v.1}$,
$e_{v.2}$ the edges containing $v$. Then by looking at the points
which span $\Delta^\vee(\T)$ we see that
$(e_{v.0}+e_{v.1}+e_{v.2})(\Delta^\vee(\T))\geq -3/2$ with the
equality only for the point $-(e^*_{v.0}+e^*_{v.1}+e_{v.2}^*)/2$ which
therefore is a vertex.  Similarly,
$(e_{v.0}+e_{v.1}-e_{v.2})(\Delta^\vee(\T))\leq 3/2$ with the equality
only for $(e^*_{v.0}+e^*_{v.1}-e_{v.2}^*)/2$.

\begin{lemma}\label{polar-polytopes}
Let $\widehat\sigma=\sum_{e\in\E} e$ then $4\Delta(\T)-2\sigma$ and
$\Delta^\vee(\T)$ are dual, or polar, one to another in the sense that
$$
\begin{array}{c}
\Delta^\vee(\T)=\left\{w\in N_\R: 
w(4\Delta(\T)-2\widehat\sigma)\geq -1\right\}\\
4\Delta(\T)-2\widehat\sigma=\left\{u\in M_\R: 
u(\Delta^\vee(\T))\geq -1\right\}
\end{array}
$$
\end{lemma}

\begin{pf}
The first equality is a restatement of \ref{defining-equations}, the
second equality follows because the polar polytope of the polar is the
original polytope, \cite[Sect 1.5]{fulton}.
\end{pf}

\begin{notation}
For a vertex of $\Delta(\T)$ we define its dual face
$u^\perp=\Delta^\vee(\T)\cap\{w: w(4u-2\widehat\sigma)=-1\}$. By
$\widetilde{u}^\perp$ we will understand the polytope which is the
convex hull of $u^\perp$ and $0\in N_\R$ while by $\widehat{u}^\perp$
we will understand the cone spanned in $N_\R$ by $u^\perp$.
\end{notation}

Let $u$ be a vertex of $\Delta(\T)$ which we can represent as a
network of paths, $\Gamma(u)$.  Then $v(u)$ is either $0$ or $2$,
depending on whether $\Gamma(u)$ contains $v$ and, similarly $e^*(u)$
is, respectively $0$ or $1$. Thus $(-v/2)(4u-2\widehat\sigma)=-1$ if
$v$ is in $\Gamma(u)$ and $(-v/2)(4u-2\widehat\sigma)=3$ otherwise. On
the other hand $(v/2-e_{v.0}^*)(4u-2\widehat\sigma)=-1$ if either $v$
is not in $\Gamma(u)$ or if both $v$ and $e^*_{v.0}$ are in
$\Gamma(u)$. Finally, $(v/2-e_{v.0}^*)(4u-2\widehat\sigma)=3$ if $v$
is in $\Gamma(u)$ but $e^*_{v.0}$ is not.

Therefore, for any vertex $u$ of $\Delta(\T)$ and any node $v\in\N$
exactly three of the following four points
$-v/2=-(e^*_{v.0}+e^*_{v.1}+e_{v.2}^*)/2$,
$v/2-e_{v.0}^*=(e_{v.1}^*+e_{v.2}^*-e_{v.0}^*)/2$,
$v/2-e_{v.1}^*=(e_{v.0}^*+e_{v.2}^*-e_{v.1}^*)/2$,
$v/2-e_{v.2}^*=(e_{v.0}^*+e_{v.1}^*-e_{v.2}^*)/2$ are in $u^\perp$
which therefore has $3n$ vertices.

\begin{ex}\label{ex-division-5-leaf}
We will visualize the points of $\hN$ on the graph of the tree in the
following way. Given a 3-valent node $v$ with edges $e_{v.0}$,
$e_{v.1}$, $e_{v.2}$, which for simplicity we denote just by numbers
on the graph, the point $-v/2$ will be denoted by the dot at the
vertex, while the point $v/2-e_{v.0}^*$ by the secant opposing the
edge $e_{v.0}$, that is

$$
\begin{array}{cccc}
\begin{xy}<8pt,0pt>:
(0,0)*={\bullet}="0" ; (-\rootthree,1)*={} **@{-},
"0" ; (\rootthree,1)*={} **@{-},
"0" ; (0,-1.8)*={} **@{-},
(-2,1)*={1},
(2,1)*={2},
(-0.5,-1.5)*={0}
\end{xy}
&
{\rm and}
&
\begin{xy}<8pt,0pt>:
(0,0)*={}="0" ; (-\rootthree,1)*={} **@{-},
(-\halfrootthree,0.5)*={}; (\halfrootthree,0.5)*={} **@{-};
"0" ; (\rootthree,1)*={} **@{-},
"0" ; (0,-1.8)*={} **@{-},
(-2,1)*={1},
(2,1)*={2},
(-0.5,-1.5)*={0}
\end{xy},
&
{\rm respectively}
\end{array}
$$

Using this notation we can put on the same picture both, the system of
paths associated to a vertex $u$ of $\Delta(\T)$ as well as the
respective points in $u^\perp$. We put only four out of eight systems
of paths from \ref{4-leaf-tree-networks} since the other ones are
obtained by renumbering of leaves.
$$
\begin{array}{cccccccc}
\begin{xy}<10pt,0pt>:
(0,0)*={}="0" ; (-1,\rootthree)*={1} **@{.},
"0" ; (-1,-\rootthree)*={2} **@{.},
"0" ; (2,0)*={}="1" **@{.},
"1" ; (3,\rootthree)*={3} **@{.},
"1" ; (3,-\rootthree)*={4} **@{.},
(1,0.5)*={0},
(-0.5,-\halfrootthree)*={}="a"; (-0.5,\halfrootthree)*={}="b" **@{-},
"a" ; (0.8,0)*={}="c" **@{-},
"b" ; "c" **@{-},
(2.5,-\halfrootthree)*={}="d"; (2.5,\halfrootthree)*={}="e" **@{-},
"d" ; (1.2,0)*={}="f" **@{-},
"e" ; "f" **@{-}
\end{xy}
&\ \ \ \ &
\begin{xy}<10pt,0pt>:
(0,0)*={\bullet}="0" ; (-1,\rootthree)*={1} **@{-},
"0" ; (-1,-\rootthree)*={2} **@{-},
"0" ; (2,0)*={}="1" **@{.},
"1" ; (3,\rootthree)*={3} **@{.},
"1" ; (3,-\rootthree)*={4} **@{.},
(1,0.5)*={0},
(-0.5,-\halfrootthree)*={}="a"; (-0.5,\halfrootthree)*={}="b" **@{},
"a" ; (0.8,0)*={}="c" **@{-},
"b" ; "c" **@{-},
(2.5,-\halfrootthree)*={}="d"; (2.5,\halfrootthree)*={}="e" **@{-},
"d" ; (1.2,0)*={}="f" **@{-},
"e" ; "f" **@{-}
\end{xy}
&\ \ \ \ &
\begin{xy}<10pt,0pt>:
(0,0)*={\bullet}="0" ; (-1,\rootthree)*={1} **@{-},
"0" ; (-1,-\rootthree)*={2} **@{-},
"0" ; (2,0)*={\bullet}="1" **@{.},
"1" ; (3,\rootthree)*={3} **@{-},
"1" ; (3,-\rootthree)*={4} **@{-},
(1,0.5)*={0},
(-0.5,-\halfrootthree)*={}="a"; (-0.5,\halfrootthree)*={}="b" **@{},
"a" ; (0.8,0)*={}="c" **@{-},
"b" ; "c" **@{-},
(2.5,-\halfrootthree)*={}="d"; (2.5,\halfrootthree)*={}="e" **@{},
"d" ; (1.2,0)*={}="f" **@{-},
"e" ; "f" **@{-}
\end{xy}
&\ \ \ \ &
\begin{xy}<10pt,0pt>:
(0,0)*={\bullet}="0" ; (-1,\rootthree)*={1} **@{-},
"0" ; (-1,-\rootthree)*={2} **@{.},
"0" ; (2,0)*={\bullet}="1" **@{-},
"1" ; (3,\rootthree)*={3} **@{-},
"1" ; (3,-\rootthree)*={4} **@{.},
(1,0.5)*={0},
(-0.5,-\halfrootthree)*={}="a"; (-0.5,\halfrootthree)*={}="b" **@{-},
"a" ; (0.8,0)*={}="c" **@{-},
"b" ; "c" **@{},
(2.5,-\halfrootthree)*={}="d"; (2.5,\halfrootthree)*={}="e" **@{-},
"d" ; (1.2,0)*={}="f" **@{-},
"e" ; "f" **@{}
\end{xy}
\end{array}
$$ 
In each of these cases the polytope $\widetilde u^\perp$ can be
divided into two simplexes, each of them having edges which make a
basis of the lattice $\hN$. For example:
$$
\begin{array}{cccccccccccccc}
\begin{xy}<10pt,0pt>:
(0,0)*={}="0" ; (-1,\rootthree)*={1} **@{.},
"0" ; (-1,-\rootthree)*={2} **@{.},
"0" ; (2,0)*={}="1" **@{.},
"1" ; (3,\rootthree)*={3} **@{.},
"1" ; (3,-\rootthree)*={4} **@{.},
(1,0.5)*={0},
(-0.5,-\halfrootthree)*={}="a"; (-0.5,\halfrootthree)*={}="b" **@{-},
"a" ; (0.8,0)*={}="c" **@{-},
"b" ; "c" **@{-},
(2.5,-\halfrootthree)*={}="d"; (2.5,\halfrootthree)*={}="e" **@{-},
"d" ; (1.2,0)*={}="f" **@{-},
"e" ; "f" **@{-}
\end{xy}
&=&
\begin{xy}<10pt,0pt>:
(0,0)*={}="0" ; (-1,\rootthree)*={1} **@{.},
"0" ; (-1,-\rootthree)*={2} **@{.},
"0" ; (2,0)*={}="1" **@{.},
"1" ; (3,\rootthree)*={3} **@{.},
"1" ; (3,-\rootthree)*={4} **@{.},
(1,0.5)*={0},
(-0.5,-\halfrootthree)*={}="a"; (-0.5,\halfrootthree)*={}="b" **@{-},
"a" ; (0.8,0)*={}="c" **@{-},
"b" ; "c" **@{-},
(2.5,-\halfrootthree)*={}="d"; (2.5,\halfrootthree)*={}="e" **@{-},
"d" ; (1.2,0)*={}="f" **@{},
"e" ; "f" **@{-}
\end{xy}
&\cup&
\begin{xy}<10pt,0pt>:
(0,0)*={}="0" ; (-1,\rootthree)*={1} **@{.},
"0" ; (-1,-\rootthree)*={2} **@{.},
"0" ; (2,0)*={}="1" **@{.},
"1" ; (3,\rootthree)*={3} **@{.},
"1" ; (3,-\rootthree)*={4} **@{.},
(1,0.5)*={0},
(-0.5,-\halfrootthree)*={}="a"; (-0.5,\halfrootthree)*={}="b" **@{-},
"a" ; (0.8,0)*={}="c" **@{-},
"b" ; "c" **@{-},
(2.5,-\halfrootthree)*={}="d"; (2.5,\halfrootthree)*={}="e" **@{-},
"d" ; (1.2,0)*={}="f" **@{-},
"e" ; "f" **@{}
\end{xy}
&=&
\begin{xy}<10pt,0pt>:
(0,0)*={}="0" ; (-1,\rootthree)*={1} **@{.},
"0" ; (-1,-\rootthree)*={2} **@{.},
"0" ; (2,0)*={}="1" **@{.},
"1" ; (3,\rootthree)*={3} **@{.},
"1" ; (3,-\rootthree)*={4} **@{.},
(1,0.5)*={0},
(-0.5,-\halfrootthree)*={}="a"; (-0.5,\halfrootthree)*={}="b" **@{-},
"a" ; (0.8,0)*={}="c" **@{},
"b" ; "c" **@{-},
(2.5,-\halfrootthree)*={}="d"; (2.5,\halfrootthree)*={}="e" **@{-},
"d" ; (1.2,0)*={}="f" **@{-},
"e" ; "f" **@{-}
\end{xy}
&\cup&
\begin{xy}<10pt,0pt>:
(0,0)*={}="0" ; (-1,\rootthree)*={1} **@{.},
"0" ; (-1,-\rootthree)*={2} **@{.},
"0" ; (2,0)*={}="1" **@{.},
"1" ; (3,\rootthree)*={3} **@{.},
"1" ; (3,-\rootthree)*={4} **@{.},
(1,0.5)*={0},
(-0.5,-\halfrootthree)*={}="a"; (-0.5,\halfrootthree)*={}="b" **@{-},
"a" ; (0.8,0)*={}="c" **@{-},
"b" ; "c" **@{},
(2.5,-\halfrootthree)*={}="d"; (2.5,\halfrootthree)*={}="e" **@{-},
"d" ; (1.2,0)*={}="f" **@{-},
"e" ; "f" **@{-}
\end{xy}
\end{array}
$$

The first equality means that $\widetilde{u}^\perp$ in this case is a union
of a simplex with edges $(e^*_1+e^*_2-e^*_0)/2$,
$(e^*_0+e^*_2-e^*_1)/2$, $(e^*_0+e^*_1-e^*_2)/2$,
$(e^*_3+e^*_4-e^*_0)/2$, $(e^*_0+e^*_3-e^*_4)/2$ and another one with
edges $(e^*_1+e^*_2-e^*_0)/2$, $(e^*_0+e^*_2-e^*_1)/2$,
$(e^*_0+e^*_1-e^*_2)/2$, $(e^*_3+e^*_4-e^*_0)/2$,
$(e^*_0+e^*_4-e^*_3)/2$.  The common part of these two simplexes is a
simplex with edges $(e^*_1+e^*_2-e^*_0)/2$, $(e^*_0+e^*_2-e^*_1)/2$,
$(e^*_0+e^*_1-e^*_2)/2$, $(e^*_3+e^*_4-e^*_0)/2$, which contains
$e^*_0/2=\left((e^*_0+e^*_3-e^*_4)/2+(e^*_0+e^*_4-e^*_3)/2\right)/2$.

This example is even more transparent when we write $\hN$ as a sum of
a rank 2 lattice spanned by $(e^*_1+e^*_2-e^*_0)/2$ and
$(e^*_3+e^*_4-e^*_0)/2$, and of rank 3 lattice spanned by
$(e^*_0+e^*_2-e^*_1)/2$, $(e^*_0+e^*_1-e^*_2)/2$,
$(e^*_0+e^*_3-e^*_4)/2$ which contains also $(e^*_0+e^*_4-e^*_3)/2$.
Then our division of the cone $\widehat{u}^\perp$ comes by multiplying
by the cone $\R_{\geq 0}(e^*_1+e^*_2-e^*_0) + \R_{\geq
0}(e^*_3+e^*_4-e^*_0)$ the standard division of the 3-dimensional cone
generated by $(e^*_0+e^*_2-e^*_1)/2$, $(e^*_0+e^*_1-e^*_2)/2$,
$(e^*_0+e^*_3-e^*_4)/2$ and $(e^*_0+e^*_4-e^*_3)/2$, see
\cite[p. 49]{fulton}, which in geometric terms is a small resolution
of a 3-dimensional quadric cone singularity giving rise to so-called
Atiyah flop.

The same argument works whenever $\Gamma(u)$ does not contain
$e_0$. Then $u^\perp$ contains $(e^*_0+e^*_2-e^*_1)/2$,
$(e^*_0+e^*_1-e^*_2)/2$, $(e^*_0+e^*_4-e^*_3)/2$,
$(e^*_0+e^*_3-e^*_4)/2$ and we can make a similar division of
$\widetilde{u}^\perp$ using the equality
$$(e^*_0+e^*_2-e^*_1)/2+(e^*_0+e^*_1-e^*_2)/2=
(e^*_0+e^*_4-e^*_3)/2+(e^*_0+e^*_3-e^*_4)/2$$

If $\Gamma(u)$ contains $e_0$ then we use the identity
$$-(e^*_0+e^*_1+e^*_2)/2+(e^*_1+e^*_2-e^*_0)/2=
-(e^*_0+e^*_3+e^*_4)/2+(e^*_3+e^*_4-e^*_0)/2$$
which presents $-e^*_0/2\in u^\perp$ as an average 
of two different pairs of vertices
to make a similar decomposition
$$
\begin{array}{cccccccccc}
\begin{xy}<10pt,0pt>:
(0,0)*={\bullet}="0" ; (-1,\rootthree)*={1} **@{-},
"0" ; (-1,-\rootthree)*={2} **@{.},
"0" ; (2,0)*={\bullet}="1" **@{-},
"1" ; (3,\rootthree)*={3} **@{-},
"1" ; (3,-\rootthree)*={4} **@{.},
(1,0.5)*={0},
(-0.5,-\halfrootthree)*={}="a"; (-0.5,\halfrootthree)*={}="b" **@{-},
"a" ; (0.8,0)*={}="c" **@{-},
"b" ; "c" **@{},
(2.5,-\halfrootthree)*={}="d"; (2.5,\halfrootthree)*={}="e" **@{-},
"d" ; (1.2,0)*={}="f" **@{-},
"e" ; "f" **@{}
\end{xy}
&=&
\begin{xy}<10pt,0pt>:
(0,0)*={\bullet}="0" ; (-1,\rootthree)*={1} **@{-},
"0" ; (-1,-\rootthree)*={2} **@{.},
"0" ; (2,0)*={}="1" **@{-},
"1" ; (3,\rootthree)*={3} **@{-},
"1" ; (3,-\rootthree)*={4} **@{.},
(1,0.5)*={0},
(-0.5,-\halfrootthree)*={}="a"; (-0.5,\halfrootthree)*={}="b" **@{-},
"a" ; (0.8,0)*={}="c" **@{-},
"b" ; "c" **@{},
(2.5,-\halfrootthree)*={}="d"; (2.5,\halfrootthree)*={}="e" **@{-},
"d" ; (1.2,0)*={}="f" **@{-},
"e" ; "f" **@{}
\end{xy}
&\cup&
\begin{xy}<10pt,0pt>:
(0,0)*={\bullet}="0" ; (-1,\rootthree)*={1} **@{-},
"0" ; (-1,-\rootthree)*={2} **@{.},
"0" ; (2,0)*={\bullet}="1" **@{-},
"1" ; (3,\rootthree)*={3} **@{-},
"1" ; (3,-\rootthree)*={4} **@{.},
(1,0.5)*={0},
(-0.5,-\halfrootthree)*={}="a"; (-0.5,\halfrootthree)*={}="b" **@{-},
"a" ; (0.8,0)*={}="c" **@{-},
"b" ; "c" **@{},
(2.5,-\halfrootthree)*={}="d"; (2.5,\halfrootthree)*={}="e" **@{},
"d" ; (1.2,0)*={}="f" **@{-},
"e" ; "f" **@{}
\end{xy}
&=&
\begin{xy}<10pt,0pt>:
(0,0)*={}="0" ; (-1,\rootthree)*={1} **@{-},
"0" ; (-1,-\rootthree)*={2} **@{.},
"0" ; (2,0)*={\bullet}="1" **@{-},
"1" ; (3,\rootthree)*={3} **@{-},
"1" ; (3,-\rootthree)*={4} **@{.},
(1,0.5)*={0},
(-0.5,-\halfrootthree)*={}="a"; (-0.5,\halfrootthree)*={}="b" **@{-},
"a" ; (0.8,0)*={}="c" **@{-},
"b" ; "c" **@{},
(2.5,-\halfrootthree)*={}="d"; (2.5,\halfrootthree)*={}="e" **@{-},
"d" ; (1.2,0)*={}="f" **@{-},
"e" ; "f" **@{}
\end{xy}
&\cup&
\begin{xy}<10pt,0pt>:
(0,0)*={\bullet}="0" ; (-1,\rootthree)*={1} **@{-},
"0" ; (-1,-\rootthree)*={2} **@{.},
"0" ; (2,0)*={\bullet}="1" **@{-},
"1" ; (3,\rootthree)*={3} **@{-},
"1" ; (3,-\rootthree)*={4} **@{.},
(1,0.5)*={0},
(-0.5,-\halfrootthree)*={}="a"; (-0.5,\halfrootthree)*={}="b" **@{},
"a" ; (0.8,0)*={}="c" **@{-},
"b" ; "c" **@{},
(2.5,-\halfrootthree)*={}="d"; (2.5,\halfrootthree)*={}="e" **@{-},
"d" ; (1.2,0)*={}="f" **@{-},
"e" ; "f" **@{}
\end{xy}
\end{array}
$$

\end{ex}

Now we shall show that the above discussion can be generalized to
the case of trees with more inner nodes.

\begin{lemma}
Suppose that $\T$ is a binary symmetric 3-valent tree with $n$ inner
nodes.  For any $u$, a vertex of $\Delta(\T)$ there exists a division
of $u^\perp$ (or, equivalently of $\widetilde{u}^\perp$) into a union
of $2^{n-1}$ (normalized) volume 1 simplexes. Equivalently, the cone
$\widetilde{u}^\perp$ can be divided into a union of simplicial cones
which are regular (i.e.~their generators form bases of $\hN$).
\end{lemma}

\begin{pf}
The construction of the division will proceed along an ascending
sequence of subtrees of $\T$, starting from an inner node of
$\T$. That is we have an ascending sequence of 3-valent trees
$$\T_1\subset \T_2 \subset \cdots \subset \T_{n-1}\subset \T_n=\T$$
where $\T_i$ has $i$ inner nodes and $\T_{i+1}$ is obtained from
$\T_i$ as a graft with a star 3-valent tree. Forgetting of edges which
are not in $\T_i$ gives a sequence of surjective maps
$M(\T)\ra\cdots\ra M(\T_i)\ra\cdots\ra M(\T_1)$ which implies a
sequences of inclusions
$\hN(\T_1)\subset\cdots\subset\hN(\T_i)\subset\cdots\subset\hN(\T_n)$.
The restriction of the networks of paths $u$ to $\T_i$ is a network on
$\T_i$ as well we will denote it by $u_i$. Clearly
$u^\perp\cap N(\T_i)_\R=u_i^\perp$.

Now we will define the division of $u^\perp_i$ inductively. The
polytope $u^\perp_1$ is just a simplex so let us assume that
$u_i^\perp=\sum \delta_i^j$ where $j=1,\dots,2^{i-1}$ and the
normalized volume of $\delta_i^j$ with respect to the lattice
$\hN_(\T_i)$ is 1. Let $v^i$ be an inner node of $\T_{i+1}$ which was a
leaf of $\T_i$, let $e^i_0$ be a petiole of $\T_i$ which become an inner
edge of $\T_{i+1}$ and let $e^i_1$ and $e^i_2$ are the two new petioles of
$\T_{i+1}$ which contain $v^i$.

Now we make argument as in \ref{ex-division-5-leaf}. If $e^i_0$ is in
$u$ then $-(e^i_0)^*/2\in u_i^\perp$ and we may assume that $e^i_1$ is
in $u$ and $e^i_2$ is not. Now from any simplex $\delta_i^j$ from the
original division of $u_i^\perp$ we produce two simplexes by adding a
new vertex at $((e_0^i)^*+(e_2^i)^*-(e_1^i)^*)/2$ and another one at
either $-((e_0^i)^*+(e_1^i)^*+(e_2^i)^*)/2$ or at
$((e_1^i)^*+(e_2^i)^*-(e_0^i)^*)/2$. Because 
$$-((e_0^i)^*+(e_1^i)^*+(e_2^i)^*)/2
+((e_1^i)^*+(e_2^i)^*-(e_0^i)^*)/2=-(e_0^i)^*$$ and $-(e_0^i)^*/2\in
u^\perp_i$ this defines a good division of $u^\perp_{i+1}$.

 If $e^i_0$ is not in $u$ then $(e^i_0)^*/2\in u_i^\perp$ and we make
 a similar construction but now we have to consider two cases: either
 none of $e_1^i$, $e^i_2$ is in $u$ or both are in $u$. At either case
 the discussion is similar to that we encountered in
 \ref{ex-division-5-leaf}.

\end{pf}

In terms of toric geometry the division process implies the following.

\begin{cor}
The affine toric variety associated to the cone $\widehat{u}^\perp$
has Gorenstein terminal singularities which admit a small resolution. 
\end{cor}

\begin{pf}
The toric singularities are Cohen-Macaulay and since all the
generators of the rays of $\widehat{u}^\perp$ lie on the hyperplane
$(4u-2\sigma)(\ \cdot\ )=-1$ the singularities in question are
Gorenstein. The division into regular simplicial cones involves adding
no extra ray so the respective resolution is small which also implies
that the original singularity is terminal.
\end{pf}

We note that the construction of the division certainly depends on the
choice of the root of the tree and changing the root gives a flop.

Let $\Sigma$ be a fan in $\hN_\R$ consisting of cones
$\widehat{u}^\perp$, where $u$ is a vertex of $\Delta(\T)$, and their
faces. In other words, $\Sigma$ contains cones spanned by the proper
faces of $\Delta^\vee(\T)$ (including the empty face, whose cone is
the zero cone).  Let us recall that equivariant line bundles on toric
varieties are in a standard way described by piecewise linear
functions on its fan, see \cite[Sect.~2.1]{oda}.  Setting
$\Lambda_{|\widehat{u}^\perp}=-u$ we define a continuous piecewise
linear function $\Lambda$ on the fan $\Sigma$ in $N_\R$ such that for
every $v\in\N$ and $e_v\in\E$ containing $v$ we have
$\Lambda(-v/2)=-1$ and $\Lambda(v/2-e_{v}^*)=0$.  The sections of the
bundle related to $\Lambda$, see \cite[Prop. 2.1]{oda}, are in
$\hM\cap\Delta(\T)$. Therefore the toric variety $X(\Sigma)$ given by
the fan $\Sigma$ can be identified with the original variety
$X(\Delta(\T))$ and the line bundle associated to $\Lambda$ is
$\cO_X(1)$. On the other hand the function $4\Lambda-2\sigma$ assumes
value 1 on the primitive vectors in rays of $\Sigma$ which allows us
to identify the canonical divisor of $X(\Delta)$, see \cite[Sect
2.1]{oda}. The result is the following.

\begin{thm}\label{Fano}
Let $\T$ be a 3-valent binary symmetric tree. Then the variety $X(\T)$
Gorenstein and Fano with terminal singularities.  Moreover it is of
index 4, that is the canonical divisor $K_{X(\T)}$ is linearly
equivalent to $\cO_{X(\T)}(-4)$.
\end{thm}

We note the following consequence of Kodaira-Kawamata-Viehweg vanishing, see
e.g.~\cite[Sect.2.5]{mori-kollar}
\begin{cor}\label{vanishing}
In the above situation $H^i(X(\T),\cO(d))=0$ for $i>0$ and $d\geq -3$.
In particular for $d\geq 0$ we have $\dim_\C H^0(X(\T),\cO(d))=h_{X(\T)}(d)$
where the latter is Poincare-Hilbert polynomial of $(X(\T),\cO(1))$.
\end{cor}


\subsection{Mutation of a tree, deformation of a model.}\label{sect-mutation-deformation}


In example \ref{4-leaf-tree-relations} we noted that a four-leaf
3-valent tree can be labeled in three non-equivalent ways. We can revert
it to say that given four numbered leaves we have three 3-valent
labeled trees connecting these leaves. By grouping in pairs the leaves
whose petioles are attached to common inner nodes we can list these as
follows: $(1,2)(3,4)$, $(1,3)(2,4)$, and $(1,4)(2,3)$. 

Now, given four pointed trees $\T_i$, where $i=1,\dots,4$ we can
produce a tree $\T$ by grafting the tree $\T_i$ along the $i$-th leaf
of a labeled 3-valent 4-leaf tree $\T_0$. Here are possible
configurations, $e_0$ denotes the inner edge of the tree $\T_0$

$$
\begin{array}{ccccc}
\begin{xy}<10pt,0pt>:
(0,0)*={}="0" ; (-0.5,\halfrootthree) **@{-}, (-1,1.5)*={\T_1},
"0" ; (-0.5,-\halfrootthree) **@{-}, (-1,-1.5)*={\T_2},
"0" ; (2,0)*={}="1" **@{-}, (1,0.5)*={e_0},
"1" ; (2.5,\halfrootthree) **@{-}, (3,1.5)*={\T_3},
"1" ; (2.5,-\halfrootthree) **@{-}, (3,-1.5)*={\T_4},
\end{xy}
&\phantom{xxxxxx}&
\begin{xy}<10pt,0pt>:
(0,0)*={}="0" ; (-0.5,\halfrootthree) **@{-}, (-1,1.5)*={\T_1},
"0" ; (-0.5,-\halfrootthree) **@{-}, (-1,-1.5)*={\T_3},
"0" ; (2,0)*={}="1" **@{-}, (1,0.5)*={e_0},
"1" ; (2.5,\halfrootthree) **@{-}, (3,1.5)*={\T_2},
"1" ; (2.5,-\halfrootthree) **@{-}, (3,-1.5)*={\T_4},
\end{xy}
&\phantom{xxxxxx}&
\begin{xy}<10pt,0pt>:
(0,0)*={}="0" ; (-0.5,\halfrootthree) **@{-}, (-1,1.5)*={\T_1},
"0" ; (-0.5,-\halfrootthree) **@{-}, (-1,-1.5)*={\T_4},
"0" ; (2,0)*={}="1" **@{-}, (1,0.5)*={e_0},
"1" ; (2.5,\halfrootthree) **@{-}, (3,1.5)*={\T_2},
"1" ; (2.5,-\halfrootthree) **@{-}, (3,-1.5)*={\T_3},
\end{xy}
\end{array}
$$

\begin{defin}\label{defin-mutation}
In the above situation we say that there exists an {\em elementary
mutation} along $e_0$ from one of the above trees to the other two. (We
note that a mutation may actually yield an equivalent tree.) We say
that two trees are {\em mutation equivalent} if there exists a
sequence of elementary mutations from one to the other.
\end{defin}

\begin{lemma}\label{trees-mutation-equivalent}
Any two 3-valent trees with the same number of leaves are mutation
equivalent.
\end{lemma}

\begin{pf}
We prove, by induction, that any 3-valent tree is mutation equivalent
to a caterpillar. To get the induction step it is enough to note that
the graft of a caterpillar tree pointed at one of its legs with a star
3 valent tree contains a distinguished inner edge the mutation of
which gives a caterpillar.
\end{pf}

Now, let us recall the basics regarding deforming subvarieties in the
projective space. Let $\cB$ be an irreducible variety (possibly
non-complete). Consider the product $\P^m\times\cB$ with the
respective projections $p_{\P}$ and $p_\cB$. Suppose that
$\cX\subset\P^m\times\cB$ is a subscheme such that the
induced projection $p_{\cB|\cX}: \cX\ra\cB$ is proper and
flat. Suppose that for two points $a, \ b\in \cB$ the respective
scheme-theoretic fibers $X_a=\cX_a$ and $X_b=\cX_b$ are reduced and
irreducible. Then we say that the subvariety $X_a$ in $\P^m$ can be
deformed to $X_b$ over the base $\cB$. This gives rise to a notion of
deformation equivalent subvarieties of $\P^m$.

\begin{defin}\label{defin-deform-equivalent}
Given two subvarieties $X_1,\ X_2$ in $\P^m$ we say that they are
deformation equivalent if their classes are in the same connected
component of the Hilbert scheme of $\P^m$.
\end{defin}

Complete intersections of the same type are deformation equivalent.
Let us consider an fundamental example, understanding of which is
essential for the proof of the main result of this section.

\begin{ex}\label{ex-def-inP7}
Let us consider $\P^7$ with homogeneous coordinates indexed by sockets
of a 4-leaf tree $\T_0$, as in example \ref{4-leaf-tree-relations}. In
$\P^7$ we consider a family of intersections of 2 quadrics
parametrized by an open subset $\cB$ of $\P^2$ with coordinates
$[t_{(12)(34)},\ t_{(13)(24)},\ t_{(14)(23)}]$. We set
$\cB=\P^2\setminus\{[1,\varepsilon,\varepsilon^2]: \varepsilon^3=1 \}$
and over $\cB$ we consider $\cX^0$ given in $\cB\times\P^7$ by
equations
$$
\begin{array}{rcl}
t_{(12)(34)}\cdot x_{1100}x_{0011} +t_{(13)(24)}\cdot
x_{1010}x_{0101}\\ +t_{(14)(23)}\cdot x_{1001}x_{0110}&
=&\left(t_{(12)(34)}+t_{(13)(24)}+t_{(14)(23)}\right)
x_{0000}x_{1111}\\ \\ 
\left(t_{(13)(24)}-t_{(14)(23)}\right)\cdot x_{1100}x_{0011}\\ 
+ \left(t_{(14)(23)}-t_{(12)(34)}\right)\cdot x_{1010}x_{0101}\\ 
+\left(t_{(12)(34)}-t_{(13)(24)}\right)\cdot
x_{1001}x_{0110} &=&0
\end{array}
$$ Three special fibers of the projection $\cX^0\ra\cB$, namely
$\cX^0_{[1,0,0]}$, $\cX^0_{[0,1,0]}$ and $\cX^0_{[0,0,1]}$, are
varieties associated to three 4-leaf trees labeled by $(12)(34)$,
$(13)(24)$ and $(14)(23)$, respectively. On the other hand $\cX^0$ is a
complete intersection of two quadrics and the map $\cX^0\ra\cB$ is
equidimensional. The latter statement follows because over $\cB$ the matrix
$$
\left[
\begin{array}{cccc}
t_{(12)(34)}&t_{(13)(24)}&t_{(14)(23)}&t_{(12)(34)}+t_{(13)(24)}+t_{(14)(23)}\\t_{(13)(24)}-t_{(14)(23)}&t_{(14)(23)}-t_{(12)(34)}&t_{(12)(34)}-t_{(13)(24)}&0
\end{array}
\right]
$$ is of rank 2 hence any fiber over $\cB$ is a complete intersection
of two non-proportional quadrics. Hence $\cX^0\ra\cB$ is flat because of
\cite[Thm. 18.16]{eisenbud}.

By $T_0\subset T_N$ denote the 4-dimensional subtorus associated to
the lattice spanned by leaves, that is a subtorus of $T_N$ with
coordinates $\chi^{v_i^*}$, where $v_i$, $i=1,\dots,4$ are leaves of
$\T_0$.  Torus $T_0$ acts on $\P^7\times\cB$ via the first coordinate,
that is, for a leaf $v_i$ of $\T_0$ and a socket $\kappa$ we have
$\lambda_{v_i}(t)(x_\kappa,t_{(.)(.)})=
t^{\kappa(v_i)}x_\kappa,t_{(.)(.)}$ Then by looking at the equations
defining $\cX^0$ wee see that the inclusion
$\cX^0\hookrightarrow\P^7\times\cB$ is equivariant with respect to
this action.

We also note that a rational map $\P^7-\ra\P^3$, regular outside
16 linear $\P^3$'s, which is given by four quadrics:
$$[x_\kappa]\ra [x_{0000}x_{1111}, x_{0011}x_{1100}, x_{0101}x_{1010},
x_{0110}x_{1001}]$$ defines a good quotient with respect to the action
of $T_0$ on $\P^7$, c.f.\cite[7.1.1]{abb}. If we take a
subvariety $\cZ^0$ in the product $\P^3\times\cB$ defined by the
equations
$$
\begin{array}{c}
t_{(12)(34)}\cdot z_1 +t_{(13)(24)}\cdot z_2+t_{(14)(23)}\cdot z_3
=(t_{(12)(23)}+t_{(13)(24)}+t_{(14)(23)})\cdot z_0\\ 
(t_{(23)(14)}-t_{(14)(23)})\cdot z_1 
+(t_{(14)(23)}-t_{(12)(34)})\cdot z_2 
+(t_{(12)(34)}-t_{(23)(14)})\cdot z_3 =0
\end{array}
$$ then $\cZ^0\ra\cB$ is equidimensional and $\cX^0$ is the fiber
product of $\P^7-\ra\P^3$ and $\cZ_0\ra\P^3$. As the result the induced
rational map $\cX^0-\ra\cZ^0$ defines a good quotient of $\cX^0$ with
respect to the action of $T_0$, \cite[7.1.4]{abb}.
\end{ex}

In what follows we construct an ambient variety which contains as
locally complete intersections a flat family of varieties containing a
geometric model of tree as well as models of the tree's elementary
mutations.

\begin{const}\label{const-quotient-of-product}
Let $\T$ be a tree with an inner edge $e_0$ which contains two
3-valent inner vertices.  We can write $\T$ as a graft of five trees:
a labeled tree $\T_0$ with four leaves $v_i$, $i=1,\dots 4$, containing
$e_0$ as an inner edge and four pointed trees $(\T_i,\ell_i)$, with
$i=1,\dots4$ which are attached to $\T_0$ along the respectively
labeled leaves. The edges in $\T$ which have common nodes with $e_0$ we
denote, respectively, by $e_i$, each $e_i$ comes from a petiole of
$\ell_i$ (or $v_i$). Recall, see \ref{polytope-of-graft-is-fiberproduct}, that
$M(\T)$ and $\Delta(\T)$ can be expressed as fiber product of
$M(\T_i)$ and $\Delta(\T_i)$, respectively. That is,
$$
\begin{array}{cc}
M(\T)=\mathop{\product}_{i=0}^4 M(\T_i)\cap \bigcap_{i=1}^4\ker(\ell_i-v_i),&
\Delta(\T)=\product_{i=0}^4 \Delta(\T_i)\cap\bigcap_{i=1}^4\ker(\ell_i-v_i)
\end{array}
$$

Now, as in \ref{toric-intersection}, we consider the lattice
$\widetilde{M}_0$ spanned on the non-trivial sockets of the tree
$\T_0$ together with the unit simplex $\widetilde{\Delta}_0\subset
\widetilde{M}_0\otimes\R$ and the maps $\widetilde{M}_0\ra M_0$ and
$\widetilde{\Delta}_0\ra \Delta_0$ which give the inclusion
$X(\T_0)\subset\P^7$ as a complete intersection of two quadrics.
Forms $v_i$, $i=1,\dots 4$ pull-back to $\widetilde{M}_0$ and we
denote them by $\tilde{v}_i$, respectively. Now we define
$$
\begin{array}{cc}
\overline{M}=\widetilde{M}_0\times\mathop{\product}_{i=1}^4 
M(\T_i)\cap \bigcap_{i=1}^4\ker(\ell_i-\tilde{v}_i)&{\rm and}\\  \\
\overline{\Delta}=\widetilde{\Delta}_0\times\product_{i=1}^4 
\Delta(\T_i)\cap\bigcap_{i=1}^4\ker(\ell_i-\tilde{v}_i)
\end{array}
$$ As in \ref{const-toric-var} we define the toric variety
$\cY=X(\overline{\Delta})$. We note that, by
\ref{fiber-product-unimodular-gen} the polytope $\overline{\Delta}$ is
normal in the lattice $\widetilde{M}_0\times\mathop{\product}_{i=1}^4
\hM(\T_i)\cap \bigcap_{i=1}^4\ker(\ell_i-\tilde{v}_i)$, which is
spanned by its vertices. Also, by the construction we have the
embeddings $X(\T)\hookrightarrow \cY\hookrightarrow \P(W^\rho_\cL)$.
\end{const}

\begin{lemma}\label{quotient-of-product}
The inclusions 
$$\begin{array}{ccc}
\overline{M}\hookrightarrow 
\widetilde{M}_0\times\mathop{\product}_{i=1}^4 M(\T_i)&
{\rm and}&
\overline{\Delta}\hookrightarrow\widetilde{\Delta}_0\times\product_{i=1}^4 
\Delta(\T_i)
\end{array}$$
induce a rational map 
$$\P^7\times \mathop{\product}_{i=1}^4 X(\T_i) - \ra \cY$$ which is a
good quotient map (of the set over which it is defined) with respect
to the action of the 4-dimensional torus $T_0$ generated by
1-parameter groups $\lambda_{v_i-\ell_i}$, where $i=1,\dots 4$. The
subvariety
$$\widehat{\cX}=\cX^0\times\mathop{\product}_{i=1}^4
X(\T_i)\hookrightarrow \cB\times\P^7\times\mathop{\product}_{i=1}^4
X(\T_i)$$ is $T_0$ equivariant and its quotient $\cX$ is
locally complete intersection in $\cB\times\cY$.
\end{lemma}
\begin{pf}
The first (quotient) part is the same as what we claim in
\ref{thm-graft-as-quotient}, this time however we repeat the argument
for all four fiber products in question.. The invariance of the
variety $\widehat{\cX}$ follows by the invariance of
$\cX^0\hookrightarrow\cB\times\P^7$ which we discussed in
\ref{ex-def-inP7}. Finally, since $\widehat{\cX}$ is a complete
intersection in $\cB\times\P^7\times\mathop{\product}_{i=1}^4 X(\T_i)$
its image $\cX$ is a locally complete intersection in the quotient
which is $\cB\times\cY$, this follows from the definition of good
quotient which locally is an affine quotient, \cite[Ch. 5]{abb},
hence functions defining $\widehat{\cX}$ locally descend to
functions defining $\cX$.
\end{pf}

\begin{lemma}\label{flatness-of-family}
Over an open set $\cB'\subset\P^2$ containing points $[1,0,0]$,
$[0,1,0]$, $[0,0,1]$ the projection morphism $\cX\ra\cB'$ is flat. The
fibers over points $[1,0,0]$, $[0,1,0]$, $[0,0,1]$ are reduced and
isomorphic to, respectively, the geometric model of $\T$ and of its
elementary mutations along the edge $e_0$.
\end{lemma}
\begin{pf}
First we note that the fibers in question, $\cX_{[*,*,*]}$, of
$\cX\ra\cB$ are geometric models as we claim. Indeed this follows from
the universal properties of good quotients, c.f.~\cite{abb}, as they
are quotients of the respective products
$\cX^0_{[*,*,*]}\times\product_{i=1}^4 X(\T_i)$, which are located, as
three invariant subvarieties, in
$\widehat\cX=\cX^0\times\product_{i=1}^4 X(\T_i)$.  This, in
particular, implies that the respective fibers of $\cX\ra\cB$ are of
the expected dimension, hence they are contained in a set
$\cB'\subset\P^2$ over which the map in question is equidimensional.
Since $\cY$ is toric it is Cohen-Macaulay and because $\cX$ is locally
complete intersection in $\cY$, it is Cohen-Macaulay too
\cite[Prop.~18.13]{eisenbud}. Finally, the map $\cX\ra\cB'$ is
equidimensional hence it is flat, because $\cB'$ is smooth, see
\cite[Thm.~18.16]{eisenbud}

\end{pf}

\begin{thm}\label{models-deformation-equivalent}
Geometric models of 3-valent trees with the same number of leaves are
deformation equivalent in $\P(W_\cL^\rho)$.
\end{thm}

\begin{pf}
This is a combination of \ref{trees-mutation-equivalent} and of
\ref{flatness-of-family}.
\end{pf}


\subsection{Hilbert-Ehrhard polynomial.}


\begin{defin}
Given two pointed trees $(\T_1,\ell_1)$ and $(\T_2,\ell_2)$ 
we define a pointed graft
which is a pointed tree
$(\T,o)=(\T_1,\ell_1)\star(\T_2,\ell_2) $ where $\T
=\T_1\phantom{*}_{\ell_1}\!\!\!\vee_{o_1} \T^{\!\!*3}
\phantom{*}_{o_2}\!\!\!\vee_{\ell_2}\T_2$, and $o$, $o_1$ and $o_2$
are the leaves of $\T^{\!\!*3}$. 
\end{defin}

\begin{ex} Pointed graft of two 3-valent stars
$$
\begin{xy}<8pt,0pt>:
(0,0)*={}="1" ; (-0.5,\halfrootthree)*={} **@{-},
"1" ; (-0.5,-\halfrootthree)*={} **@{-},
"1" ; (0.8,0) **@{-},(1,0)*={\circ} 
\end{xy}\ 
\star\ 
\begin{xy}<8pt,0pt>:
(0,0)*={}="2" ; (0.5,\halfrootthree)*={} **@{-},
"2" ; (0.5,-\halfrootthree)*={} **@{-},
"2" ; (-0.8,0) **@{-},(-1,0)*={\circ}
\end{xy}\ 
=\ 
\begin{xy}<8pt,0pt>:
(0,-0.5)*={}="0" ; (-0.9,0)*={}="1"  **@{-},
"1" ; (-1.4,\halfrootthree)*={} **@{-},
"1" ; (-1.4,-\halfrootthree)*={} **@{-},
"0" ; (0.9,0)*={}="2" **@{-},
"2" ; (1.4,\halfrootthree)*={} **@{-},
"2" ; (1.4,-\halfrootthree)*={} **@{-},
"0" ; (0,-1.2) **@{-}, (0,-1.4)*={\circ}
\end{xy}
$$
\end{ex}

By arguments used in the proof of
\ref{polytope-of-graft-is-fiberproduct} we also get.

\begin{prop}\label{polytope-of-graft-is-simpleproduct2}
Let $(\T_1,\ell_1)$ and $(\T_2,\ell_2)$ be two pointed trees. Then

$$\Delta(\T_1\phantom{*}_{\ell_1}\!\!\star_{\ell_2}\T_2)=
\Delta(\T_1)\phantom{*}_{\ell_1}\!\!\!\times_{o_1}
\Delta(\T^{\!\!*3})\phantom{*}_{o_2}\!\!\!\times_{\ell_2}\Delta(\T_2)$$
\end{prop}

Let us consider a 3-dimensional lattice $M=\Z e_0\oplus\Z e_1\oplus\Z
e_2$ with a fixed tetrahedron $\Delta^0$ with vertices $0,\ e_0+e_1,\
e_0+e_2,\ e_1+e_2$. By $\widehat{M}\subset M$ we denote the index 2
sublattice spanned on the vertices of $\Delta^0$. 

\begin{defin}
Let $n$ be a positive integer and let $f_1^n=f_1,\ f_2^n=f_2$ be two
functions defined on the set $\{0,\dots n\}$ with values in $\Z$ or,
more generally, in an arbitrary ring or algebra (we use the
superscript $^n$ to indicate the domain of $f$'s).  For any
$k\in\{0,\dots,n\}$ we define
$$
(f_1\star f_2)(k)=\sum_{\stackrel{u\in\widehat{M}\cap n\Delta^0}{e_0^*(u)=k}}\ 
\left( f_1(e_1^*(u))\cdot f_2(e_2^*(u))\right) 
$$
\end{defin}

We note that $\star$ is commutative, that is $f^n_1\star
f^n_2=f^n_2\star f^n_1$, but possibly not associative. By
$(f^n)^{\star m}$ we denote the $\star$ product of $m$ copies of a
chain of $f^n$, that is $f^n\star(f^n\star(\dots (f^n\star
f^n)\dots)$.  By $1^n$ we denote the constant function
$\{0,\dots,n\}\ra\{1\}\subset\Z$.

 A function $f^n: \{0,\dots, n\}\ra\Z$ will be called {\em symmetric}
if $f^n(k)=f^n(n-k)$.

\begin{lemma}\label{counting-points}
If $f_1=f^n_1,\ f_2=f^n_2: \{0,\dots, n\}\ra\Z$ are symmetric
functions then $f_1\star f_2$ is a symmetric function as well and
moreover for $k\leq n/2$ we have
$$
(f_1\star f_2)(k)=
2\cdot\left(\sum_{i=0}^{k-1}\sum_{j=0}^{i} f_1(i)f_2(k+i-2j)\right)
+\left(\sum_{i=k}^{n-k}\sum_{j=0}^k f_1(i)f_2(k+i-2j)\right)
$$
In particular, for $k\leq n/2$
$$
(f_1*1)(k)=2\sum_{i=0}^{k-1}(i+1)f_1(i)+\sum_{i=k}^{n-k}(k+1)f_1(i)
$$
\end{lemma}

\begin{pf}

Let us look at the sections of the tetrahedron $n\Delta^0$ with
hyperplanes $(e_0^*)^{-1}(k)$. We picture the situation for $n=6$ and
$k=0,\dots 6$, the dotted square is the section of the cube with the
lower left corner satisfying relation $e_1^*=e^*_2=0$, the section of
the tetrahedron denoted with solid line and points inside the
(closed) tetrahedron denoted by $\bullet$.
$$
\begin{array}{ccccccc}
k=0&k=1&k=2&k=3&k=4&k=5&k=6\\
\begin{xy}<8pt,0pt>:
(0,0)*={\bullet}="0" ; (0,6)*={}="1" **@{.}, 
"0" ; (6,0)*={}="2" **@{.},
"2" ; (6,6)*={\bullet}="3" **@{.},
"1" ; "3" **@{.},
"0" ; "3" **@{-},
(1,1)*={\bullet}, (2,2)*={\bullet},(3,3)*={\bullet},(4,4)*={\bullet},(5,5)*={\bullet}
\end{xy}
&
\begin{xy}<8pt,0pt>:
(0,0)*={}="0" ; (0,6)*={}="1" **@{.}, 
"0" ; (6,0)*={}="2" **@{.},
"2" ; (6,6)*={}="3" **@{.},
"1" ; "3" **@{.},
(1,0)*={\bullet}="4" ; (0,1)*={\bullet}="5" **@{-},
"4" ; (6,5)*={\bullet}="6" **@{-},
"5" ; (5,6)*={\bullet}="7" **@{-},
"6" ; "7" **@{-},
(2,1)*={\bullet}, (3,2)*={\bullet}, (4,3)*={\bullet}, 
"7" ; (5,4)*={\bullet} **@{.},
"4" ; (1,2)*={\bullet} **@{.}, 
(2,3)*={\bullet}, (3,4)*={\bullet}, (4,5)*={\bullet}
\end{xy}
&
\begin{xy}<8pt,0pt>:
(0,0)*={}="0" ; (0,6)*={}="1" **@{.}, 
"0" ; (6,0)*={}="2" **@{.},
"2" ; (6,6)*={}="3" **@{.},
"1" ; "3" **@{.},
(2,0)*={\bullet}="4" ; (0,2)*={\bullet}="5" **@{-},
"4" ; (6,4)*={\bullet}="6" **@{-},
"5" ; (4,6)*={\bullet}="7" **@{-},
"6" ; "7" **@{-},
(3,1)*={\bullet}, 
"7" ; (4,2)*={\bullet} **@{.}, 
(5,3)*={\bullet}, 
(1,3)*={\bullet}, 
"4" ; (2,4)*={\bullet} **@{.}, 
(3,5)*={\bullet}, 
(1,1)*={\bullet}, (2,2)*={\bullet},(3,3)*={\bullet},(4,4)*={\bullet},(5,5)*={\bullet}
\end{xy}
&
\begin{xy}<8pt,0pt>:
(0,0)*={}="0" ; (0,6)*={}="1" **@{.}, 
"0" ; (6,0)*={}="2" **@{.},
"2" ; (6,6)*={}="3" **@{.},
"1" ; "3" **@{.},
(3,0)*={\bullet}="4" ; (0,3)*={\bullet}="5" **@{-},
"4" ; (6,3)*={\bullet}="6" **@{-},
"5" ; (3,6)*={\bullet}="7" **@{-},
"6" ; "7" **@{-},
(4,1)*={\bullet}, (5,2)*={\bullet}, 
(1,4)*={\bullet}, (2,5)*={\bullet}, 
(1,2)*={\bullet}, (2,1)*={\bullet},
(4,5)*={\bullet},(5,4)*={\bullet},
(2,3)*={\bullet}, (3,2)*={\bullet},
(3,4)*={\bullet},(4,3)*={\bullet},
\end{xy}
&
\begin{xy}<8pt,0pt>:
(0,0)*={}="0" ; (0,6)*={}="1" **@{.}, 
"0" ; (6,0)*={}="2" **@{.},
"2" ; (6,6)*={}="3" **@{.},
"1" ; "3" **@{.},
(4,0)*={\bullet}="4" ; (0,4)*={\bullet}="5" **@{-},
"4" ; (6,2)*={\bullet}="6" **@{-},
"5" ; (2,6)*={\bullet}="7" **@{-},
"6" ; "7" **@{-},
(3,1)*={\bullet}, (2,2)*={\bullet}, (1,3)*={\bullet}, 
(5,3)*={\bullet}, (4,4)*={\bullet}, (3,5)*={\bullet}, 
(5,1)*={\bullet}, (4,2)*={\bullet},(3,3)*={\bullet},(2,4)*={\bullet},(1,5)*={\bullet}
\end{xy}
&
\begin{xy}<8pt,0pt>:
(0,0)*={}="0" ; (0,6)*={}="1" **@{.}, 
"0" ; (6,0)*={}="2" **@{.},
"2" ; (6,6)*={}="3" **@{.},
"1" ; "3" **@{.},
(5,0)*={\bullet}="4" ; (0,5)*={\bullet}="5" **@{-},
"4" ; (6,1)*={\bullet}="6" **@{-},
"5" ; (1,6)*={\bullet}="7" **@{-},
"6" ; "7" **@{-},
(4,1)*={\bullet}, (3,2)*={\bullet}, (2,3)*={\bullet}, (1,4)*={\bullet},
(5,2)*={\bullet}, (4,3)*={\bullet}, (3,4)*={\bullet}, (2,5)*={\bullet}
\end{xy}
&
\begin{xy}<8pt,0pt>:
(0,0)*={}="0" ; (0,6)*={\bullet}="1" **@{.}, 
"0" ; (6,0)*={\bullet}="2" **@{.},
"2" ; (6,6)*={}="3" **@{.},
"1" ; "3" **@{.},
"1" ; "2" **@{-},
(5,1)*={\bullet}, (4,2)*={\bullet},(3,3)*={\bullet},(2,4)*={\bullet},(1,5)*={\bullet}
\end{xy}
\end{array}
$$ The definition of $f_1^n\star f_2^n$ is sum of the product of
$f^n_i$'s over the lattice points of such a section. The sections over
$k$ and $n-k$ are obtained by a reflection with respect to either
$e^*_1=1/2$ or $e^*_2=1/2$. Thus if one of $f^n_i$'s is symmetric then
the $f^n_1\star f^n_2$ is symmetric as well.

On the other hand for $0\leq k\leq n-k$ the tetrahedron's section is a
rectangle with vertices $(k,0),\ (0,k),\ (n-k,n),\ (n,n-k)$ which we
divide into two triangles and a parallelogram, the division is
indicated by dotted vertical line segments for boxes labeled by
$k=1,2$ in the above diagram. Because functions $f^n_i$ are symmetric
the values of the product $f^n_1\cdot f^n_2$ are the same for the
points which are central symmetric with respect to the center of the
square.  Thus in the formula of the lemma we take the value
$f_1(a)f_2(b)$ for all integral pairs $(a,b)$ in the left hand side
triangle and multiply it by 2 (that is the first summand in the
formula) and add the sum over the parallelogram.
\end{pf}

\begin{ex}\label{ex-counting-points}
We note that $(1^n)^{\star 2}(k)=(k+1)(n-k+1)$ is the number of
lattice points in the rectangle used in the argument in the above
proof of \ref{counting-points}. On the other hand by using the formula
from \ref{counting-points} one gets
$$
(1^n)^{\star 3}(k)=
{1\over 6}(k+1)(n-k+1)(n^2+kn-k^2+5n+6)
$$\end{ex}

Let us recall that given a lattice polytope $\Delta\subset M_\R$ for
any positive integer $n$ we define Ehrhard function $h_\Delta$
as follows: 
$$
h_\Delta(n)=|\left(n\cdot \Delta\cap M\right)|
$$ If $\Delta$ satisfies the assumptions of \ref{const-toric-var} then
$h_\Delta=h_{X(\Delta)}$ where the latter is the Poincare-Hilbert
polynomial of $(X(\Delta),\cO(1))$ which, by definition, is equal to
$\dim_\C H^0(X(\Delta),\cO(m))$ for $m\gg 0$.

\begin{defin}\label{def-relative-ehrhard}
Let $\Delta\subset M_\R$ be a lattice polytope which is not contained
in any hyperplane and let $v\in N$ be a non-zero form on $M$. Suppose
that $v(\Delta)\subset [0,1]$. We define its relative Ehrhard function
$f^n_{\Delta,v}: \{0,\dots n\}\ra\Z$ by setting
$$
f^n_{\Delta,v}(k)=|v^{-1}(k)\cap n\cdot\Delta\cap M|
$$
\end{defin}

We note that, clearly, $\sum_{k=0}^n f_{\Delta,v}^n(k)=h_\Delta(n)$ is
the usual Ehrhard function.  Thus, in case of \ref{const-toric-var}
the above definition can be restated in purely geometric fashion.
\begin{lemma}
Suppose that $\Delta$ satisfies assumptions of \ref{const-toric-var}
and $v$ is as in \ref{def-relative-ehrhard}.  Let us consider a
linearization of the action of the 1-parameter group $\lambda_v$ on
$H^0(X(\Delta),\cO(n))$ which has non-negative weights and the
eigenspace of the zero weight is nontrivial. Then $f^n_{\Delta,v}(k)$
is equal to the dimension of the eigenspace of the action of
$\lambda_v$ of weight $k$.
\end{lemma}

\begin{pf}
This is a consequence of the standard properties of $X(\Delta)$,
\ref{toric-via-polytope}.4.
\end{pf}

\begin{lemma}
Let $(\T_1,\ell_1)$ and $(\T_2,\ell_2)$ be two pointed trees and let
$f^n_{\ell_1}$ and $f^n_{\ell_2}$ be two relative Ehrhard functions
associated to $\Delta(\T_1)$ and $\Delta(\T_2)$, respectively.  
If $(\T,o)=(\T_1,\ell_1)\star(\T_2,\ell_2)$ and $f^n_o$ is the relative
Ehrhard function associated to $\Delta(\T)$ then
$f^n_o=f^n_{\ell_1}\star f^n_{\ell_2}$
\end{lemma}

\begin{pf}
The definitions of $\star$ are made accordingly.
\end{pf}

\begin{ex}\label{easy-ehrhard-poly}
By using \ref{ex-counting-points} we find out that
$$\sum_{k=0}^n(1^n)^{\star 2}(k)={{(n+1)(n+2)(n+3)}\over{6}}$$
which is the Poincare-Hilbert polynomial of $(\P^3,\cO(1))$
while 
$$\sum_{k=0}^n(1^n)^{\star 3}(k)=
{{(n+1)(n+2)(n+3)(n^2+4n+5)}\over{30}}$$ which is Poincare-Hilbert
polynomial of intersection of two quadrics in $\P^7$.
\end{ex}

\begin{thm}\label{thm-associativity}
Let us consider three pointed trees $(\T_i,\ell_i)$, with $i=1,2,3$
with relative Ehrhard functions $f^n_i=f^n_{\ell_i}$ associated to
polytopes $\Delta(\T_i)$, respectively. Then
$$
\left(f^n_{1}\star f^n_{2}\right)\star f^n_{3}=
f^n_{1}\star \left(f^n_{2}\star f^n_{3}\right)
$$
\end{thm}

\begin{pf}
Let $\ell$ denote the distinguished leaf of the result of the $\star$
operation on the trees. Then the relative Ehrhard function
$\left(f^n_{1}\star f^n_{2}\right)\star f^n_{3}$ and, respectively,
$f^n_{1}\star \left(f^n_{2}\star f^n_{3}\right)$ is related to one of
the following trees, each of them is obtained by an elementary mutation
from the other:
$$
\begin{array}{ccc}
\begin{xy}<10pt,0pt>:
(0,0)*={}="0" ; (-0.5,\halfrootthree) **@{-}, (-1,1.5)*={\ell},
"0" ; (-0.5,-\halfrootthree) **@{-}, (-1,-1.5)*={\T_3},
"0" ; (1,0)*={}="1" **@{-},
"1" ; (1.5,\halfrootthree) **@{-}, (2,1.5)*={\T_1},
"1" ; (1.5,-\halfrootthree) **@{-}, (2,-1.5)*={\T_2},
\end{xy}
&\longleftrightarrow&
\begin{xy}<10pt,0pt>:
(0,0)*={}="0" ; (-0.5,\halfrootthree) **@{-}, (-1,1.5)*={\ell},
"0" ; (-0.5,-\halfrootthree) **@{-}, (-1,-1.5)*={\T_1},
"0" ; (1,0)*={}="1" **@{-},
"1" ; (1.5,\halfrootthree) **@{-}, (2,1.5)*={\T_3},
"1" ; (1.5,-\halfrootthree) **@{-}, (2,-1.5)*={\T_2},
\end{xy}
\end{array}
$$ 

Now we repeat the construction \ref{const-quotient-of-product},
with obvious modifications. Namely, we define 
a polytope 
$$\overline{\Delta}=\widetilde{\Delta}_0\times\product_{i=1}^3
\Delta(\T_i)\cap\bigcap_{i=1}^3\ker(\ell_i-\tilde{v}_i)$$ where
$\widetilde{\Delta}_0$ is the unit simplex as in
\ref{toric-intersection}. We define a toric variety
$\cY=X(\overline{\Delta})$ with the embedding in $\P(W^\rho_\E)$ and
the action of the group $\lambda_\ell$. 

Next, as in \ref{quotient-of-product} we define a subvariety
$\cX\subset \cB\times\cY$ such that the projection $p_{\cB}:
\cX\ra\cB$ is flat and its two fibers are varieties associated to the
above two pointed trees, see \ref{flatness-of-family}.  Because of the
flatness the sheaf $(p_\cB)_*(p_\cY^*(\cO(n))$ is locally free for
each $n\geq 0$, see \cite[III.9.9, III.12.9]{hartshorne} and
\ref{vanishing}. Moreover, by the construction, the action of the
group $\lambda_\ell$ on $\cY$ leaves $\cX\subset\cB\times\cY$, as we
noted in \ref{ex-def-inP7}.  Finally, the decomposition into
eigenspaces of the action of $\lambda_\ell$ on $H^0(\cY,\cO(n))$
restricts into a respective eigenspace decomposition of the action of
$\lambda_\ell$ on fibers of $(p_\cB)_*(p_\cY^*(\cO(n))$, which are
equal to $H^0(\cX_b,\cO(n))$, for $b\in\cB$. This implies that the
dimension of the respective eigenspaces is locally constant,with
respect to the parameter $b\in\cB$ hence the relative Ehrhard function
of fibers of $p_\cB$ is constant which concludes the argument.
\end{pf}

Let us underline the fact that although the invariance of the Hilbert
polynomial is a standard property of a flat family the above result is
about the invariance of the family with respect to an action of a
1-parameter group, the group $\lambda_\ell$ in our case. 

The above theorem \ref{thm-associativity} implies that the operation
$\star$ on relative Ehrhard functions of polytopes of 3-valent trees
is not only commutative (which is obvious from its definition) but
also associative. This implies that the function does not depend on
either the shape nor the location of the leaf. More precisely we have
the following formula which allows to compute the Hilbert-Ehrhard
polynomial very efficiently.

\begin{cor}\label{relative-ehrhard-equal}
If $(\T,\ell)$ is a pointed 3-valent tree with $r+1$ leaves then
$$f^n_{\Delta(\T),\ell}=(1^n)^{\star r}$$ 
\end{cor}

\renewcommand{\thesection}{A} 

\section{Appendix}

\subsection{Normal polytopes, unimodular covers}\label{sect-unimodular}

A lattice simplex $\Delta^0\subset {M}_\R$ with vertices
$v_0,\dots v_r$ is called unimodular if vectors
$v_1-v_0,\dots,v_r-v_0$ span $M$. We say at a lattice polytope
$\Delta\subset {M}_\R$ has a unimodular covering if
$\Delta=\bigcup_\nu\Delta^0_\nu$ where $\Delta^0_\nu$ are unimodular
simplexes. This definition is taken from \cite{bgt} where we also have
the following result.

\begin{lemma}\label{unimodular=>normal}
If a lattice polytope $\Delta\subset {M}_\R$ has a unimodular
covering then it is normal.
\end{lemma}

The following observation is probably known but we include its proof
because of the proof of the subsequent lemma.

\begin{lemma}\label{product-unimodular}
Let $\Delta_1\subset (M_1)_\R$ and $\Delta_2\subset (M_2)_\R$ be two
unimodular simplexes.  Then $\Delta_1\times\Delta_2$ has a
unimodular covering in $M_1\times M_2$.
\end{lemma}

\begin{pf}
We can assume that $\Delta_1$ has vertices 
$0, e_1,\dots,e_r$ and $\Delta_2$ has vertices
$0, f_1,\dots,f_s$. Suppose that $x\in (M_1)_\R\times(M_2)_\R$
is as follows:
$$
x=\sum_{i=1}^r a_ie_i + \sum_{i=j}^s b_jf_j
$$
where $a_i, b_j\geq 0$ and $\sum a_i \leq \sum b_j\leq 1$. 

The union of unimodular simplexes contained in
$\Delta_1\times\Delta_2$ is a closed subset. Therefore if $x$ is
not contained in any modular subsimplex of
$\Delta_1\times\Delta_2$ then any small perturbation of $x$ has
this property as well. Thus we are free to assume that all $a_i$'s and
$b_j$'s are nonzero and any two non-empty subsets of $a_i$'s and
$b_j$'s have different sum, in particular $a_{1}+\dots+a_{p}\ne
b_{1}+\dots+b_{q}$ for any reasonable $(p,q)$.  Let $m$ be such
$b_1+\dots+b_{m-1}<a_1+\dots+a_r<b_1+\dots+b_{m}$. We set
$b'_{m}=(b_1+\dots+b_{m})-(a_1+\dots+a_r)$.

In order to prove the lemma we will find $r+m-1$ positive numbers
$c_{i,j}$ indexed by some pairs
$(i,j)\in\{1,\dots,r\}\times\{1,\dots,m\}$ such that 
$$x=\sum_{ij} c_{i,j}(e_i+f_j)+
b'_{m}f_{m}+b_{m+1}f_{m}+\dots+b_sf_s$$ and the set of respective
vectors $(e_i+f_j)$ together with $f_m,\dots f_s$ can be modified via
addition or subtraction of pairs among them to the standard basis
$e_1,\dots e_r, f_1,\dots, f_s$.

The coefficients $c_{i,j}$ are defined inductively according to the
following rules. The first coefficient is
$c_{1,1}=\min\{a_1,b_1\}$. Suppose that the last defined coefficient is
$c_{i_0,j_0}$. If $(i_0,j_0)=(r,m)$ then we are done so assume that it
is not the case. Then, because of our assumption that the sequences $(a_i)$
and $(b_j)$ have no equal partial sums, either 
 $a_1+\dots+a_{i_0} > b_1+\dots +b_{j_0}$, or 
 $a_1+\dots+a_{i_0} < b_1+\dots +b_{j_0}$.
In the former case we set 
$$c_{i_0,j_0+1}=\min\{b_{j_0+1},\ (a_1+\dots+a_{i_0})-(b_1+\dots+b_{j_0})\}$$
whereas in the latter case we define
$$c_{i_0+1,j_0}=\min\{a_{i_0+1},\
(b_1+\dots+b_{j_0})-(a_1+\dots+a_{i_0})\}$$ 

The verification that $\sum_{i=1}^r c_{i,j}=b_j$ for $j=1\dots m-1$
and $\sum_{j=1}^m c_{i,j}=a_i$ for $i=1,\dots r$ is easy and left for
the reader.  Similarly, a simple backtracking allows to modify the set
of the respective vectors $e_i+f_j$ with $f_m,\dots f_s$ to the
standard basis for $M_1\times M_2$.
\end{pf}

\begin{lemma}\label{fiber-product-unimodular}
Let $\Delta_1\subset (M_1)_\R$, $\Delta_2\subset(M_2)_\R$ be two
unimodular simplexes. We consider two homomorphisms $\ell_i: M_i\ra\Z$
such that $(\ell_i)_\R(\Delta_i)\subset [0,1]$. The the fiber product
$\Delta=(\Delta_1)_{\ell_1}\!\!\times\!_{\ell_2}(\Delta_2)$ has a
unimodular covering with respect to the fiber product lattice
$M=(M_1)_{\ell_1}\!\!\times\!_{\ell_2}(M_2)$.
\end{lemma}

\begin{pf}
The argument is a variation of the one used in the previous lemma. We
can assume that $\Delta_1$ has vertices $0, e^0_1,\dots,e^0_{r_0},
e^1_1,\dots,e^1_{r_1}$ and $\Delta_2$ has vertices $0,
f^0_1,\dots,f^0_{s_0},f^1_1,\dots,f^1_{s_1}$ where
$\ell_1(e^0_i)=\ell_2(f^0_j)=0$ and $\ell_1(e^1_i)=\ell_2(f^1_j)=1$
for suitable $i$'s and $j$'s.  Suppose that $x\in
(M_1)_\R\times(M_2)_\R$ is as follows:
$$ 
x=\sum_{i=1}^{r_0} a^0_ie^0_i+ \sum_{i=1}^{r_1} a^1_ie^1_i +
\sum_{j=1}^{s_0} b^0_jf^0_j+ \sum_{j=1}^{s_1} b^1_jf^1_j
$$ where $a^0_i, a^1_i, b^0_j, b^1_j\geq 0$, $\sum a^0_i +\sum
a^1_i\leq 1$, $\sum b^0_j+\sum b^1_j\leq 1$ and moreover $\sum
a^1_i=\sum b^1_j $.  The latter conditions ensures that
$\ell_1(x)=\ell_2(x)$ and it is the only condition which can not be
made perturbed, as in the proof of the previous lemma.

We write $x=x_0+x_1$ where $x_0=\sum a^0_ie^0_i+ \sum b^0_jf^0_j$ and
$x_1= \sum a^1_ie^1_i+\sum b^1_jf^1_j$ and we repeat the proof of
\ref{product-unimodular} for $x_0$ and $x_1$ separately.  The only
difference is that, because of the equality $\sum a^1_i=\sum b^1_j$,
the construction will give $r_1+s_1-1$ coefficients $c^1_{i,j}$ and
associated pairs of vectors $e^1_i+f^1_j$ which will enable to write
$x_1=\sum c^1_{i,j}(e^1_i+f^1_j)$. Thus, clearly, the respective
vectors $e^1_i+f^1_j$ do not constitute a basis of the lattice spanned
by $e^1_1,\dots,e^1_{r_1},f^1_1,\dots, f^1_{s_1}$ but of this lattice
intersected with $\ker(\ell_1-\ell_2)$. That is, among the chosen
$r_1+s_1-1$ vectors $e^1_i+f^1_j$ we have $e^1_1+f^1_1$ and
$e^1_{r_1}+f^1_{s_1}$ and if $e^1_i+f^1_j$ is among them then either
$e^1_{i+1}+f^1_j$ or $e^1_i+f^1_{j+1}$ is among them as well (but not
both). We are to prove that any $e^1_i+f^1_j$ can be obtained as a sum
of them.  But this follows because
$$\left(e^1_{i}+f^1_{j}\right)+\left(e^1_{i+1}+f^1_{j+1}\right)=
\left(e^1_{i+1}+f^1_{j}\right)+\left(e^1_{i}+f^1_{j+1}\right)$$ so any one
of the above above four vectors is a combination of the other three
and this observation can be used repeteadly to complete our claim.
\end{pf}

\begin{cor}\label{fiber-product-unimodular-gen}
Let $\Delta_1\subset (M_1)_\R$, $\Delta_2\subset(M_2)_\R$ be two
polytopes which have covering by unimodular simplexes. We consider two
homomorphisms $\ell_i: M_i\ra\Z$ such that
$(\ell_i)_\R(\Delta_i)\subset [0,1]$. Then the fiber product
$\Delta=(\Delta_1)_{\ell_1}\!\!\times\!_{\ell_2}(\Delta_2)$ has a
unimodular covering with respect to the fiber product lattice
$M=(M_1)_{\ell_1}\!\!\times\!_{\ell_2}(M_2)$
\end{cor}

\begin{pf}
The fiber product of $\Delta_1$ and $\Delta_2$ is covered by fiber
products of simplexes from the unimodular cover of each of them. Thus
the result follows by \ref{fiber-product-unimodular}.
\end{pf}

Since the polytope of the star 3-valent tree is a unit tetrahedron we
get the following.

\begin{prop}\label{tree-polytope-normal}
If $\T$ is a binary symmetric {\bf 3-valent} tree then its polytope in
$\widehat{M}(\T)$ has unimodular covering hence it is normal.
\end{prop}


\subsection{Two 3-valent trees with 6 leaves}

One of the fundamental questions regarding the phylogenetic trees is
the following. Given two (3-valent binary symmetric) trees $\T_1$ and
$\T_2$ suppose that $\Delta(\T_1)\iso\Delta(\T_2)$ as lattice
polytopes, or the projective models $X(\T_1)$ and $X(\T_2)$ are
projectively equivalent. Does it imply that the trees are equivalent
(as CW complexes) as well?

We tackled the problem by understanding the difference of models of
the two simplest non-equivalent trees. These are 6-leaf trees pictured
below, respectively, a 3-caterpillar tree and a tree which we
call a snow flake, \cite{ss}.
$$
\begin{xy}<15pt,0pt>:
(0,0)*={}="0" ; (-\halfroottwo,\halfroottwo)*={} **@{-},
"0" ; (-\halfroottwo,-\halfroottwo)*={} **@{-},
"0" ; (1,-0.1)*={}="1" **@{-},
"1" ; (0.9,-1)*={}     **@{-},
"1" ; (2,0.1)*={}="2" **@{-},
"2" ; (2.1,-1)*={}     **@{-},
"2" ; (3,0)*={}="5"  **@{-},
"5" ; (\threeandhalfroottwo,\halfroottwo)*={} **@{-},
"5" ; (\threeandhalfroottwo,-\halfroottwo)*={} **@{-},
\end{xy}
\phantom{xxxxxxxxxxx}
\begin{xy}<12pt,0pt>:
(0,-0.2)*={}="0" ; (-0.9,0.1)*={}="1"  **@{-},
"1" ; (-1.4,\halfrootthree)*={} **@{-},
"1" ; (-1.4,-\halfrootthree)*={} **@{-},
"0" ; (0.9,0.1)*={}="2" **@{-},
"2" ; (1.4,\halfrootthree)*={} **@{-},
"2" ; (1.4,-\halfrootthree)*={} **@{-},
"0" ; (0,-1.2)*={}="3" **@{-},
"3" ; (-0.7,-1.8) **@{-},
"3" ; (+0.7,-1.8) **@{-} 
\end{xy}
$$ The snow flake tree is obtained from the 3-caterpillar tree by
elementary mutation along its middle inner edge. Therefore their
Hilbert-Ehrhard polynomials are equal and computed with [{\tt maxima}]
to be as follows.
$$
\begin{array}{rcl}
h(n)&=&{1\over22680}\left(n+1\right)\left(n+2\right)\left(n+3\right)\cdot\\
&&\left(31\,n^6+372\,n^5+1942\,n^4+5616\,n^3+9511\,n^2+8988\,n+3780\right)
\end{array}
$$
On the other hand we can distinguish their polytopes in terms of some
combinatorial invariants.

\begin{ex} 
Given a polytope $\Delta$ we define its incidence matrix $(a_{ij})$ as
follows: $(a_{ij})$ is a symmetric matrix with integral entries such
that for $i\leq j$ the number $a_{ij}$ is equal to the number of
$i$-dimensional faces contained in $j$-dimensional faces of
$\Delta$. In particular $a_{ii}$ is the number of $i$-dimensional
faces. The following is the incidence matrix of a polytope of the
snow flake tree.
$$
\begin{array}{ccccccccc}
32&  480&   2400&   6144&   9312&   8832&   5280&   1920&   384\\
480& 240&  2400&   9456& {\bf19920}&{\bf24960}&{\bf19200}&{\bf8880}&{\bf2256}\\
2400& 2400& 760& 5944&{\bf19008}&{\bf32552}&{\bf32408}&{\bf18792}&{\bf5872}\\
6144& 9456& 5944& 1316&{\bf8400}&{\bf21744}&{\bf29308}&{\bf21720}&{\bf8388}\\
9312&{\bf19920}&{\bf19008}&{\bf8400}& 1392&{\bf7200}&{\bf14640}&{\bf14640}&{\bf7200}\\
8832&{\bf24960}&{\bf32552}&{\bf21744}&{\bf7200}& 940&{\bf3820}&{\bf5760}&{\bf3820}\\
5280&{\bf19200}&{\bf32408}&{\bf29308}&{\bf14640}&{\bf3820}& 406& 1224& 1224\\
1920&{\bf8880}&{\bf18792}&{\bf21720}&{\bf14640}&{\bf5760}& 1224& 108& 216\\
384&{\bf2256}&{\bf5872}&{\bf8388}&{\bf7200}&{\bf3820}& 1224& 216& 16
\end{array}
$$
And this is the incidence matrix of the polytope of a 3-caterpillar tree.
$$
\begin{array}{ccccccccc}
32& 480& 2400& 6144& 9312& 8832& 5280& 1920& 384\\
480& 240& 2400& 9456&{\bf19904}&{\bf24896}&{\bf19104}&{\bf8816}&{\bf2240}\\
2400& 2400& 760& 5944&{\bf18976}&{\bf32408}&{\bf32168}&{\bf18616}&{\bf5824}\\
6144& 9456& 5944& 1316&{\bf8384}&{\bf21648}&{\bf29112}&{\bf21552}&{\bf8336}\\
9312&{\bf19904}&{\bf18976}&{\bf8384}& 1392&{\bf7184}&{\bf14584}&{\bf14576}&{\bf7176}\\
8832&{\bf24896}&{\bf32408}&{\bf21648}&{\bf7184}& 940&{\bf3816}&{\bf5752}&{\bf3816}\\
5280&{\bf19104}&{\bf32168}&{\bf29112}&{\bf14584}&{\bf3816}& 406& 1224& 1224\\
1920&{\bf8816}&{\bf18616}&{\bf21552}&{\bf14576}&{\bf5752}& 1224& 108& 216\\
384&{\bf2240}&{\bf5824}&{\bf8336}&{\bf7176}&{\bf3816}& 1224& 216& 16
\end{array}
$$ 

Both matrices were computed by [{\tt polymake}].  We note that
although both polytopes have the same number of faces of respective
dimension their incidences are different (indicated in boldface).

\end{ex}


\subsection{Volume distribution}

The leading coefficient in the Ehrhard polynomial of a lattice
polytope $\Delta$ can be identified as the volume of $\Delta$ (with
respect to the lattice in question, whose unit cube is assumed to have
volume 1). Similarly, we can define a relative volume function which
will measure the distribution of the volume of $\Delta(\T)$ with
respect to a leaf $\ell$ of $\T$.  Because of
\ref{relative-ehrhard-equal} this function does not depend either on
the shape of the tree nor on the choice of the leaf $\ell$. Moreover
we will normalize it so that its integral over the unit segment is 1.

\begin{figure}[h]
\center{\includegraphics[height=8cm]{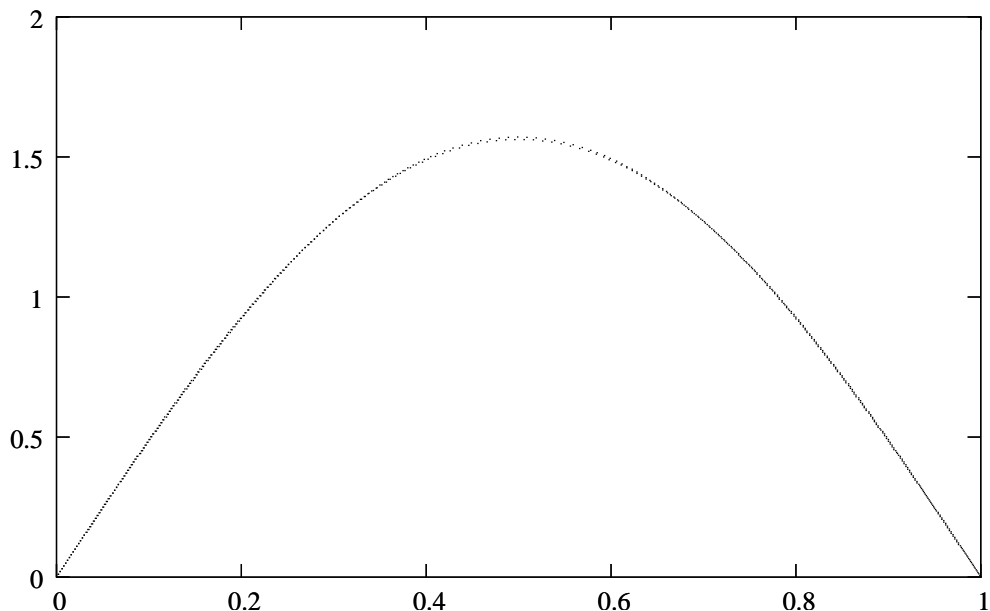}}
\caption{Polynomials $\delta^2$ and $\delta^{100}$ at the same
diagram, by \cite{gnuplot}.}
\end{figure}

If $\delta^r: [0,1]\ra\R$ is the normalized volume distribution with
respect to a leaf of a 3-valent tree with $r$ leaves then because of
\ref{counting-points} we have $\delta^n(t)=\delta^n(1-t)$ and for
$t\in(0,1/2)$ we get the following recursive formula
$$ \delta^{r+1}(t)=d_{n+1}\cdot\left(
2\cdot\int_{0}^t s\cdot\delta^r(s)ds+\int_{t}^{1-t}t\cdot\delta^r(s)ds
\right)
$$ where $d_{n+1}$ is a constant such that
$\int_0^1\delta^{r+1}(s)ds=1$.  From this it follows that $\delta^r$
is a polynomial of degree $2r$. However, the numerical experiments
which we have made seem to indicate that for $r>3$ the actual values
of $\delta^r$ do not depend too much on $r$, see Fig.~1.  It seems
that this function does not see the shape of the tree (which is
because it comes from the relative Hilbert-Ehrhard polynomial) but
also almost disregards its size (or dimension of the model)



\bigskip

\noindent Authors' address: Instytut Matematyki UW, Banacha 2, 02-097
Warszawa, Poland
\noindent 
{\tt wkrych@mimuw.edu.pl}\hfill{\tt jarekw@mimuw.edu.pl}
\end{document}